\newcommand\Modified{1.2; September 24, 2010}

\documentclass[11pt]{amsart}

\usepackage{amsfonts, amsmath,amscd, amssymb, latexsym, mathrsfs,  slashed, stmaryrd, verbatim, wasysym }
\usepackage[all]{xy}
\usepackage{graphicx}
\usepackage{hyperref}

\newcommand\datver[1]{\def\datverp%
{\par\boxed{\boxed{\text{Version: #1; Run: \today}}}}}

\newtheorem*{acknowledgements}{Acknowledgements}
\newtheorem{theorem}{Theorem}
\newtheorem*{Remark}{Remark}
\newtheorem{corollary}[theorem]{Corollary}
\newtheorem{definition}[theorem]{Definition}

\newtheorem{lemma}[theorem]{Lemma}
\newtheorem{proposition}[theorem]{Proposition}
\newtheorem{remark}[theorem]{Remark}

\numberwithin{equation}{section}
\numberwithin{theorem}{section}

\newcommand\tomega{\widetilde{\omega}}
\newcommand\bomega{\bar\omega}
\mathchardef\mhyphen="2D
\newcommand\Fhc{\operatorname{F-hc}}
\newcommand\bbf{\operatorname{bf}}
\newcommand\tf{\operatorname{tf}}
\newcommand\cuf{\operatorname{cf}}

\newcommand{\df}[1]{\mathfrak{#1}}

\newcommand{\bhs}[1]{\mathfrak B_{#1}}
\renewcommand{\tilde}{\widetilde}
\renewcommand{\bar}{\overline}
\renewcommand{\Re}{\operatorname{Re}}
\renewcommand{\hat}[1]{\widehat{#1}}

\newcommand{\wt}[1]{\widetilde{#1}}
\newcommand{\rest}[1]{\big\rvert_{#1}} 

\newcommand{\Rint}{\sideset{^R}{}{\int}}

\newcommand{\ddt}[1]{\frac{\partial}{\partial #1}}
\newcommand{\dzero}[1]{\frac{\partial}{\partial #1}\biggr\rvert_{#1 =0}}
\newcommand\pa{\partial}

\newcommand\eg{e\@.g\@. }

\newcommand\eps\varepsilon

\newcommand\hc{\operatorname{hc}}

\newcommand\CI{{\mathcal{C}}^{\infty}}

\newcommand{\lrpar}[1]{\left( #1 \right)}
\newcommand{\lrspar}[1]{\left[ #1 \right]}

\newcommand{\lrbrac}[1]{\left\lbrace #1 \right\rbrace}
\newcommand{\norm}[1]{\lVert #1 \rVert}
\newcommand{\abs}[1]{\left\lvert #1 \right\rvert}

\newcommand\Area{\operatorname{Area}}

\renewcommand\det{\operatorname{det}}

\newcommand\diag{\operatorname{diag}}

\newcommand\dvol{\operatorname{dA}}

\DeclareMathOperator*{\FP}{\operatorname{FP}}
\newcommand\fun{\operatorname{F}}

\renewcommand\Re{\operatorname{Re}}
\DeclareMathOperator*\Res{\operatorname{Res}}

\newcommand\RTr[1]{{}^R\!\operatorname{Tr}\left( #1 \right)}

\newcommand\Spec{\operatorname{Spec}}

\newcommand\Tr{\operatorname{Tr}}

\newcommand\Mand{\text{ and }}
\newcommand\Mas{\text{ as }}

\newcommand\Mwith{\text{ with }}

\newcommand\paperintro%
        {%
         }
\newcommand\paperbody%
        {%
         }

\newcommand\bbH{\mathbb{H}}

\newcommand\bbM{\mathbb{M}}
\newcommand\bbN{\mathbb{N}}

\newcommand\bbR{\mathbb{R}}
\newcommand\bbS{\mathbb{S}}

\newcommand\cC{\mathcal{C}}

\newcommand\cH{\mathcal{H}}

\newcommand\cK{\mathcal{K}}
\newcommand\cL{\mathcal{L}}

\newcommand\cO{\mathcal{O}}
\newcommand\cP{\mathcal{P}}

\newcommand\cR{\mathcal{R}}
\newcommand\cS{\mathcal{S}}

\newcommand\sC{\mathscr{C}}

\newcommand\sF{\mathscr{F}}

\newcommand{\FunA}{{
\declareslashed{}{\scriptstyle{I}}{.5}{0}{\scriptstyle{I}}
\slashed{\scriptstyle{I}} }}
\newcommand{\FunB}{{
\declareslashed{}{\smile}{.45}{.65}{\FunA}
\slashed{\FunA} }}
\newcommand{\Fun}{{
\declareslashed{}{\frown}{.05}{-1.35}{\FunB}
\slashed{\FunB} }}

\DeclareMathAlphabet{\mathpzc}{OT1}{pzc}{m}{it}

\datver{\Modified}
\begin{document}

\title[Ricci flow and the determinant of the Laplacian]{Ricci flow and the determinant of the Laplacian on non-compact surfaces}
\author{Pierre Albin}
\author{Clara L. Aldana}
\author{Fr\'ed\'eric Rochon}
\address{Courant Institute of Mathematics and Institute for Advanced Study\newline
Current address: University of Illinois at Urbana-Champaign}
\email{palbin@illinois.edu}
\address{Department of Mathematics, Universidad de los Andes\newline
Current address: Max Planck Institute for Gravitational Physics}
\email{clara.aldana@aei.mpg.de}
\address{Department of Mathematics,  Australian National University\newline
Current address: Universit\'e du Qu\'ebec \`a Montr\'eal}
\email{rochon.frederic@uqam.ca}
\thanks{This material is based upon work supported by the National Science Foundation under grant DMS-0635607002 and an NSF postdoctoral fellowship (first author) and a NSERC discovery grant (third author). Any opinions, findings, and conclusions or recommendations expressed in this material are those of the authors and do not necessarily reflect the views of the NSF}


\newcommand\Sec[1]{\cS_{ #1 }} 
\newcommand\Curv[1]{\hat \cR_{ #1 }} 

\begin{abstract}
On compact
surfaces with or without boundary, Osgood, Phillips and Sarnak
proved that the maximum of the determinant of the Laplacian within a conformal class of metrics
with fixed area
occurs at a metric of constant curvature and, for negative Euler
characteristic, exhibited a flow from a given metric to a constant
curvature metric along which the determinant increases. The aim of
this paper is to perform a similar analysis for the determinant of
the Laplacian on a non-compact surface whose ends are asymptotic to
hyperbolic funnels or cusps. In that context, we show that the Ricci flow converges to a
metric of constant curvature and that the determinant increases
along this flow.
\end{abstract}

\maketitle

\tableofcontents

\paperintro
\section*{Introduction}

A non-compact surface of finite topology $M$ has a natural compactification $\bar M$ to a surface with boundary obtained by attaching a circle `at infinity' to each end.
To discuss asymptotic expansions of functions or sections of a vector bundle on $M,$ it is very convenient to pass to the compactification and make use of a {\bf boundary defining function} or {\bf bdf}.
A boundary defining function for $Y,$ a component of the boundary of $M,$ is a smooth non-negative function on $\bar M$ that is equal to zero precisely on $Y,$ and has non-vanishing differential on $Y.$
A bdf for $\pa M$ is known as a {\bf total boundary defining function}.

\begin{definition}
Let $\bar M$ be a compact surface with boundary. Assume that the connected components of the boundary of $M$ have been partitioned into `funnel ends' and `cusp ends'
\begin{equation*}
    \pa_{\fun} \bar M = \{ Y_1, \ldots, Y_{n_{\fun}} \}, \quad
    \pa_{\hc} \bar M = \{ Y_{n_{\fun}+1}, \ldots, Y_{n_{\fun}+n_{\hc}} \}.
\end{equation*}
A metric $g$ on the interior of $M$ is a {\bf funnel-cusp metric} or {\bf $\Fhc$ metric} if,
for each $Y_{i} \in \pa_{\fun}M$ there is a bdf $x_i$ and a collar neighborhood $(0,\epsilon)_{x_{i}} \times \bbS^1$ of $Y$ in $M$ on which the pull-back of $g$ is equal to
\begin{equation}\label{ModelFun}
    e^{\varphi}\lrpar{ \frac{dx_{i}^2 + d\theta_{i}^2}{x_{i}^2} },
\end{equation}
and similarly
for each $Y_{i} \in \pa_{\hc} M$ there is a bdf $x_{i}$ and a collar neighborhood $(0,\epsilon)_{x_{i}} \times \bbS^1$ of $Y$ in $M$ on which the pull-back of $g$ is equal to
\begin{equation}\label{ModelCusp}
    e^{\varphi}\lrpar{ \frac{dx_{i}^2}{x_{i}^2} + x_{i}^2 d\theta_{i}^2 }.
\end{equation}
In both cases, $d\theta_{i}^2$ is the round metric on the circle of length one, and $\varphi$ is required to be a smooth function on $\bar M$ equal to a constant at $x_{i}=0.$
\label{se.3a}\end{definition}

A simple example of $\Fhc$ metric is  the \textbf{horn},
\begin{equation}
   \sF = \bbS^1 \times (0,\infty)_s,
   \text{ with the metric } g_\sF = \frac{ds^2 + d\theta^2}{s^2}.
\label{horn.1}\end{equation}
It has a funnel end for $\{ s < 1 \}$  and a cusp end for $\{ s > 1 \}.$
In fact, the hyperbolic metric on the quotient of $\bbH^2$ by a geometrically finite discrete group of hyperbolic isometries is a $\Fhc$ metric \cite[Example 2.1]{Perry}.

The interior of any manifold with boundary can be endowed with a $\Fhc$ metric.
A result of Mazzeo and Taylor \cite[\S 2]{MazzeoTaylor}, and also a consequence of the present manuscript, is that any such metric can be conformally transformed to a hyperbolic metric.
Note that funnel ends are also referred to as conformally compact or asymptotically hyperbolic \cite{Mazzeo:Hodge}.

For each $i,$ we will assume that $x_{i}$ is equal to $1$ outside a collar
neighborhood of $Y_{i}.$  As a global boundary defining function, we can
therefore consider
\begin{equation}
   x= \prod_{i=1}^{n_{\hc}+n_{\fun}} x_{i}.
\label{bdf.1}\end{equation}
We can also take
\begin{equation}
x_{\fun}= \prod_{i=1}^{n_{\fun}}x_{i} \quad \mbox{and} \quad
x_{\hc}= \prod_{i=n_{\fun}+1}^{n_{\fun}+n_{\hc}}x_{i}
\label{bdf.2}\end{equation}
as boundary defining functions for $\pa_{\fun}\bar M$ and $\pa_{\hc}\bar M$
respectively.
Any non-compact hyperbolic metric is an example of a $\Fhc$-metric, as is any compact perturbation of a hyperbolic metric.

In analogy with the case of a manifold with boundary,
we say that the funnel ends of a $\Fhc$ metric $(M,g)$ are {\bf totally geodesic} if, in the description \eqref{ModelFun}, we have
\begin{equation}\label{TotGeodFun}
    \varphi - \varphi\rest{x_i=0} = o(x_i)
\end{equation}
and, if this condition holds at the cusp ends, we say that they are totally geodesic.  Unless otherwise stated, we will usually assume that the metrics considered are totally geodesic.

In this paper, we propose to define and study the determinant of the (positive definite) Laplacian on such surfaces.
On compact manifolds, the determinant of the Laplacian is a global spectral invariant
originally defined by Ray and Singer \cite{Ray-Singer} using the zeta
function of the Laplacian to regularize the product of its non-zero
eigenvalues. For a non-compact manifold, the Laplacian typically has both point spectrum and  continuous spectrum, which means further
regularizations are  necessary to define its determinant. One approach,
due to M\"uller \cite{Mueller:RD} and extended to our situation by
the second author \cite{Aldana} and by Borthwick, Judge, and Perry \cite{BJP:Duke}, is to define a relative determinant
by comparing the Laplacian to a model operator along the ends of the
manifold. For non-compact surfaces with constant curvature, other
approaches involve techniques from hyperbolic geometry, e.g.,
\cite{Efrat1},\cite{Efrat2} \cite{TZ1}, \cite{BJP:Selberg}. We follow a method originally due to Richard
Melrose \cite{APSBook} that allows us to use refined information
about the Schwartz kernel of the heat operator (from \cite{Albin2} and
\cite{Vaillant}) to define its {\em renormalized} trace and extend
the definition of the zeta function and determinant. This method has the advantage of being
very flexible and systematic.  It
has been used by Andrew Hassell in his proof of the Cheeger-M\"uller
theorem \cite{Hassell}, and in recent work of two of the authors
\cite{Albin-Rochon1, Albin-Rochon2, Albin-Rochon3}.

As in the compact case, our determinant admits a Polyakov formula des\-cribing the variation of
the determinant under conformal deformations
of the metric.
If all of the ends of $(M,g)$ are asymptotically cusps, then the
area of $M$ is finite and the formula is very similar to the one for the compact case.
This was done in terms of relative determinants by the second author in \cite{Aldana} and by Borthwick, Judge and Perry \cite{BJP:Duke} (using the
Mazzeo-Taylor uniformization \cite{MazzeoTaylor} as a starting point).  In
general, if there are any funnel ends, then the area of $M$ is infinite and the behavior of the determinant of the Laplacian is more complicated.
The same renormalization process used to extend the determinant of the Laplacian can be used to define renormalized integrals and, in particular, a `renormalized area', ${}^R\!\!\Area.$  Complications arise because this renormalized integral is not a positive functional. Indeed, we show below that compactly supported conformal changes can be used to make the renormalized area equal to any given real number, positive or negative!

Explicitly, our Polyakov formula (Theorem~\ref{thm:Polyakov} below) stipulates that if  $g_\tau = e^{\omega(\tau)}g_0$ is a family of $\Fhc$ metrics on a non-compact surface $M$ with totally geodesic ends, then the  determinant of the Laplacian satisfies
\begin{equation}\label{PolyakovFiniteIntro}
    \pa_{\tau} \log \det \Delta_\tau
    =
    -\frac1{24\pi} \int_M \omega'(\tau) R_\tau \dvol_{\tau}
    + \pa_\tau \log \Area_\tau(M)
\end{equation}
in the finite area case and
\begin{equation}\label{PolyakovInfiniteIntro}
    \pa_{\tau} \log \det \Delta_\tau
    =
    -\frac1{24\pi} \sideset{^R}{_M}\int  \omega'(\tau) R_\tau \dvol_{\tau}
\end{equation}
in the infinite area case, where $R_\tau$ is the scalar curvature of $g_\tau.$

By analogy with the compact case and the result of \cite{OPS2}, one would expect that among
all metrics in a given conformal class and  with fixed  (renormalized) area, the determinant should be maximal on the one with constant scalar curvature.  For the relative determinant on Riemann surfaces with cusps ends, such a result was recently obtained by the second author \cite{Aldana}. For $\Fhc$ metrics with only funnel ends and of constant curvature outside a compact set, the corresponding result was obtained by Borthwick, Judge, and Perry \cite{BJP:Duke}, again for the relative determinant.

For the determinant considered here, our strategy to establish the analog of the result of \cite{OPS2} is to use the Ricci flow.  The relevance of Ricci flow for the determinant of the Laplacian on closed surfaces was pointed out already in \cite{OPS2} by Osgood, Phillips, and Sarnak. Later M\"uller and Wendland \cite{MullerWendland}  (see also the work of Kokotov and Korotkin \cite{KokotovKorotkin}) verified that the determinant on closed surfaces increases along this flow.

On compact Riemann surfaces, Ricci flow was considered by Hamilton (see also the work of Cao \cite{Cao} for the K\"ahler-Ricci flow) who showed in \cite{Hamilton}
that on surfaces of negative Euler characteristic, the normalized Ricci flow exists for all time and converges to a metric of constant curvature.  The proof of this result elegantly follows from the study of the evolution equation of an accessory potential function.  Following the same strategy,
this result was  recently generalized by Ji, Mazzeo, and \v{S}e\v{s}um \cite{Ji-Mazzeo-Sesum} to non-compact surfaces with asymptotically cusp ends, the key new difficulty in this case being the construction
of the potential function, which turns out to be much more delicate due to the presence of cusps.  See also the work of Chau \cite{Chau} for the generalization of the result of Cao \cite{Cao} to non-compact K\"ahler manifolds.

In section~\ref{sec:Ricci}, we generalize further the result of \cite{Ji-Mazzeo-Sesum} to include also funnels (see Theorem~\ref{thm:RicciFlow}).  Along the way, we also carefully study how the asymptotic behavior of the metric evolves along the flow, an important point for the study of the determinant. A new feature in our case is that, as long as there is at least one end asymptotic to a funnel, the Euler characteristic need not be negative, but can also equal $0$ or $1.$  This is consistent with the fact that there exists a metric with negative constant scalar curvature in these cases, \eg the horn and the hyperbolic plane $\bbH^{2}.$  Again, the key step is the construction of the potential function.  Essentially by a doubling
construction along the funnels, we can reduce to a situation where there are only cusps and obtain our potential function from the one of \cite{Ji-Mazzeo-Sesum}.  Since  the area is infinite in this setting,  we have to proceed differently to define a normalized Ricci flow.  Instead of using the average scalar curvature, we can use any fixed  constant $\sC$ and consider the normalized Ricci flow
\begin{equation}
      \frac{\pa g}{\pa t}(t)= (\sC-R_{g(t)})g(t).
\label{int.1}\end{equation}
The constant $\sC$ has to be negative to insure the flow is aiming towards a metric of constant negative scalar curvature.
A particularly interesting choice is to take this number to be the renormalized average curvature, in which case the renormalized area is preserved along the flow.

Finally, in the last section, we combine our Polyakov formula with the convergence result for the Ricci flow
to get our main result, which says (see Theorem~\ref{ExtraAssumption} for the precise statement) that among all $\Fhc$ metrics g  in a given conformal class with totally geodesics ends satisfying
\[
      {}^R\!\!\Area(g) = -2\pi \chi(M)
\]
with scalar curvature asymptotically equal to $-2$ in each funnel end, the determinant of the Laplacian
is greatest at the hyperbolic metric in this conformal class.

Indeed, if $g_{0}$ is such a metric and $g(t)=e^{\omega(t)}g_{0}$ is the solution to the normalized
Ricci flow with $\sC$ given by the (renormalized) average curvature, then the variation of the determinant along the flow is nonnegative  and given by
\begin{equation}
\begin{split}
    \pa_t \log \det(\Delta_{g(t)})
    &= -\frac1{24\pi} \Rint \omega'(t) R_t \dvol_{t} \\
    &= \frac1{24\pi} \Rint (R_t - \sC)^2 \dvol_{t} + \frac{\sC}{12\pi} \Rint (R_t - \sC) \dvol_t \\
    &= \frac1{24\pi} \Rint (R_t - \sC)^2 \dvol_{t} \ge 0,
\end{split}
\label{int.2}\end{equation}
the last integral being nonnegative since it does not need to be renormalized because
$R_{t}-\sC=\cO(x_{\fun})$ along the flow (See section~\ref{max.0} for all the details).

\begin{Remark}
Related results have been obtained since the first appearance of this paper. We mention in particular the works of Gregor and Topping \cite{GregorTopping}, of Bahuaud \cite{Bahuaud} and of the third author and Zhang \cite{RochonZhang}.
\end{Remark}

\begin{acknowledgements}
The authors are happy to acknowledge the hospitality and support of MSRI through the program \emph{Analysis on Singular Spaces}, where this work was begun. We are also grateful to Rafe Mazzeo for many illuminating discussions on the Ricci flow, to David Borthwick for helpful comments, and to the anonymous referees.
\end{acknowledgements}

\paperbody
\section{Renormalization on non-compact Riemann surfaces}\label{renor.0}

Our approach to extend the definition of the determinant of the
Laplacian to non-compact surfaces is through renormalized integrals.
In this section we briefly review how to renormalize the integral of
a density on a manifold with boundary provided it has asymptotic
expansions at each boundary face. We refer the reader to
\cite{APSBook} for more details, where renormalized integrals are called
`$b$-integrals', as well as
\cite{Hassell-Mazzeo-Melrose1}, \cite{Albin1},
\cite[Appendix]{Albin-Rochon3}.

\subsection{The choice of the boundary defining function}\label{cbdf.0}

The choice of coordinates $(x_{i},\theta_{i})$ in \eqref{ModelFun}
or \eqref{ModelCusp} and in particular the choice of boundary
defining function \eqref{bdf.1} is not fixed by the conformal
structure of the metric.  A local conformal change of coordinates in
\eqref{ModelFun} or \eqref{ModelCusp} would induce a different
boundary defining function.  However, the boundary compactification
$\bar M,$ which can be seen as being obtained from $(M,g)$ via our
choice of boundary defining function \eqref{bdf.1}, does only depend
on the conformal class of $g.$

\begin{proposition}
The boundary compactification $\bar M$ does not depend on the choice of conformal coordinates
in \eqref{ModelFun} and \eqref{ModelCusp}.  If $(\hat{x}_{i},\hat{\theta}_{i})$ corresponds to a different
choice of conformal coordinates, then $\hat{x}_{i}= x_{i} h(x_{i},\theta_{i})$ for some smooth function $h$
with $h(0,\theta_{i})$ nowhere zero on the boundary component $Y_{i}.$
\label{cbdf.1}\end{proposition}
\begin{proof}
The horn $\sF$ of \eqref{horn.1} is conformal to the punctured unit disk using the complex coordinate
$\zeta= e^{2\pi i z}$ with $z=\theta+is,$ the cusp end corresponding to the puncture and the funnel end corresponding to the boundary of the unit disk.    In this coordinate, the boundary defining function of the
cusp end is given by $\rho_{\hc}= \frac{-2\pi}{\log |\zeta |}$ and by $\rho_{\fun}= -\frac{\log |\zeta|}{2\pi}$ for the funnel end.

Using this coordinate near a cusp end, we see that a local change of conformal coordinates near the cusp is given by a holomorphic function $f(\zeta)=\zeta g(\zeta)$ with
$g(0)\ne 0$ and induces the new boundary defining function
\begin{equation}
      \hat{\rho}_{\hc}= \frac{-2\pi}{\log |f(\zeta)|}.
\label{cbdf.2}\end{equation}
A quick inspection shows that $f$ extends to give a smooth function on the boundary compactification
at the cusp (defined by $\rho_{\hc}$) and that $\hat{\rho}_{\hc}=\rho_{\hc} h(\rho_{\hc},\theta)$ with
$h$ a smooth function with $h(0,\theta)$ nowhere zero.

Similarly, a local change of conformal coordinates near a funnel end is given by a holomorphic
function $f$ defined for $1-\epsilon<|\zeta| <1$ for some $\epsilon>0$ and which extends to a continuous
function on the unit circle in such a way that  $|f(\zeta)|=1$ whenever $|\zeta|=1.$  It defines a new
boundary defining function $\hat{\rho}_{\fun}= -\frac{\log |f(\zeta)|}{2\pi}.$  By the Schwarz reflection
principle, $f$ extends to be smooth on the boundary compactification (defined by $\rho_{\fun}$) and
$\hat{\rho}_{\fun}= \rho_{\fun}h(\rho_{\fun}, \theta)$ for some smooth function $h$ with
$h(0,\theta)$ nowhere zero.

\end{proof}

In this paper, we will assume that the boundary defining function \eqref{bdf.1} and the conformal
coordinates $(x_{i},\theta_{i})$ in \eqref{ModelFun} and \eqref{ModelCusp} are given and fixed.  Notice that whether or not an end is \emph{totally geodesic} \eqref{TotGeodFun} depends on this choice.
In the funnel case, this notion is  related with
the important r\^ole  played by `special' or `geodesic' bdf's (see
\cite{Graham}, \cite{Albin1}). These are bdf's $x$ such that
$|\frac{dx}x|_g$ is constant in a neighborhood of the funnel
boundary (the value {\em at} the funnel boundary is independent of
the choice of bdf). In order to make use of these results, we point
out that a bdf satisfying \eqref{TotGeodFun} is `close to'
geodesic for the associated $\Fhc$ metric. By \cite[Lemma
2.1]{Graham} we know that there exists a geodesic bdf $\hat x$ for
$g,$ unique in a neighborhood of $\pa_{\fun}M,$ such that $\hat x^2
g\rest{T \pa_{\fun}M}$ is the metric $ d\theta_{i}^{2}$ on each
component $Y_{i}$ of $\pa_{\fun}\bar M.$

\begin{lemma}\label{ConformalBdf}
Let $g$ be a $\Fhc$ metric, $x$ be a bdf satisfying the assumptions above, and let $\varphi_0 = \varphi\rest{x=0}.$
Let $\hat x$ be a geodesic bdf as described above that coincides with $x$ along the cusp ends.
If we write $\hat x = e^{\phi}x$ then $\phi = -\frac12\varphi_0 + \cO(x)$ as $x \to 0.$
If the funnel ends are totally geodesic, then $\phi = -\frac12\varphi_0 + o(x).$
\end{lemma}

\begin{proof}
Assume without loss of generality that we only have funnel ends.
We need to have $\hat x^2 g\rest{T\pa M} = e^{-\varphi_0}x^2g \rest{T\pa M},$ but since
\begin{equation*}
    \hat x^2g = e^{2\phi} x^2g
\end{equation*}
this means we must have $2\phi+\varphi_0 = \cO(x).$
We are also asking that $|\frac{d\hat x}{\hat x}|^2_{g}$ be constant near $\hat x =0.$ The constant is necessarily equal to
\begin{equation*}
    \left. \abs{\frac{dx}x}^2_g \right|_{x=0} = e^{-\varphi_0}
\end{equation*}
so we need to have
\begin{equation*}
\begin{gathered}
    e^{-\varphi_0}
    = g \lrpar{ \frac{d\hat x}{\hat x}, \frac{d\hat x}{\hat x} }
    = g \lrpar{ \frac{dx}x + d\phi, \frac{dx}x + d\phi } \\
    = \lrpar{ e^{-\varphi} + 2g\lrpar{ \frac{dx}x, d\phi} + g\lrpar{ d\phi, d\phi } }.
\end{gathered}
\end{equation*}
Since $g(d\phi, d\phi) = \cO(x^2),$ this implies $\pa_x \phi = o(1)$ if $\varphi - \varphi_0 = o(x).$
\end{proof}

Since the volume form of a $\Fhc$ metric blows-up at the funnel ends to second order, this lemma implies as we will see that, for metrics totally geodesic along the funnel ends, renormalization results that require a geodesic bdf also hold for a bdf satisfying the assumptions above.

\subsection{Renormalized integrals} $ $\newline
Let $M$ be a manifold with boundary. A function $f$ on $M$ is {\bf
polyhomogeneous} if it is smooth in the interior of $M$ and, at each
connected component $N$ of $\pa M,$ $f$ has an asymptotic expansion
in terms of a bdf $x$ for $N$ and $\log x,$
\begin{equation*}
    f \sim \sum a_{s,p} x^s (\log x)^p.
\end{equation*}
We require that, for each $\ell \in \bbN,$ there are only finitely many $a_{s,p}$ with $\Re s<\ell.$

If $\mu$ is a smooth non-vanishing density on $\bar M,$ then the asymptotic expansion of $f$ can be used to meromorphically continue the function
\begin{equation*}
    z\mapsto \int_M x^z f \mu
\end{equation*}
and then we define the renormalized integral of $f$ to be
\begin{equation*}
    \sideset{^R}{_M}\int f \mu = \FP_{z=0}\int_M x^z f \mu.
\end{equation*}
Alternately, one can use the asymptotic expansion of $f$ to show that the function
\begin{equation*}
    \eps \mapsto \int_{x \geq \eps} f \mu
\end{equation*}
has an asymptotic expansion as $\eps \to 0$ and then define
\begin{equation*}
    \sideset{^H}{_M}\int f\mu = \FP_{\eps=0}\int_{x\geq \eps} f \mu.
\end{equation*}
This method of renormalizing is often known as {\em Hadamard renormalization} and is used in \cite{APSBook} while the previous method is known as {\em Riesz renormalization} and is used in \cite{Melrose-Nistor}. In \cite[\S2.3]{Albin1} it is shown that, under certain natural conditions which will always hold in this manuscript, these two renormalizations coincide.

\subsection{The renormalized area} \label{sec:RenArea} $ $\newline
Let $(M,g)$ be a surface of infinite area.
Given a total boundary defining function on $M,$ we can define the renormalized area of $M,$ ${}^R\Area(g),$ by taking the renormalized integral of the volume form of $g.$
In sharp contrast to the usual area, the renormalized area {\em need not be positive.}

The renormalized area generally depends on the choice of bdf used to define it though, as we are working with a fixed bdf,  ${}^R\Area(g)$ is unambiguously defined for a $\Fhc$ metric. The sign of ${}^R\Area(g)$ does not reveal any information about the behavior of the metric at infinity. Indeed, the renormalized area is clearly additive under compact perturbations and one can compactly change the metric (even conformally) and arrange for the renormalized area to be any real number, positive or negative.

\subsection{The renormalized Gauss-Bonnet theorem} $ $\newline
In this section we describe the extension of the Gauss-Bonnet theorem to $\Fhc$-metrics.

Recall that the Gauss-Bonnet theorem for a surface with boundary $(M, \bar g)$ says that
\begin{equation*}
    \int_M R(\bar g) \; \dvol_{\bar g} + \int_{\pa M} \Fun = 4\pi\chi(M)
\end{equation*}
where $\Fun$ is a density on $\pa M$ built up from the second fundamental form of $\pa M$ in $M.$
In particular, if the boundary of $M$ is totally geodesic, then $\int_M R(\bar g) \; \dvol_{\bar g} = 4\pi\chi(M).$
A straight-forward computation using Lemma \ref{ConformalBdf} (see \cite[Theorem 4.5]{Albin1}) shows that, if all of the funnel ends of $g$ are totally geodesic \eqref{TotGeodFun}, and $x$ is as in \eqref{bdf.1}, then there is a renormalized Gauss-Bonnet theorem with no contribution from the boundary.

\begin{theorem}[Gauss-Bonnet for $\Fhc$ metrics]
Let $(M,g)$ be a non-compact surface with a $F,\hc$ metric. If
$\Area(M)$ is finite, then
\begin{equation*}
    \int_M R(g) \; \dvol_g = 4\pi\chi(M).
\end{equation*}
If $\Area(M)$ is infinite and the infinite volume ends of $M$ are totally geodesic, then
\begin{equation*}
    \Rint R(g) \; \dvol_g  = 4\pi\chi(M).
\end{equation*}
\end{theorem}

\begin{remark}
For exactly hyperbolic metrics this appears in \cite[\S 2]{BJP:Selberg} for surfaces and in \cite{Epstein} for higher dimensional hyperbolic manifolds without cusps; see also \cite{Albin1} for higher dimensional asymptotically hyperbolic metrics that are asymptotically Einstein.

If the funnel ends of $g$ are not totally geodesic, one can check that the boundary contribution in the corresponding renormalized Gauss-Bonnet theorem is essentially the same as the boundary contribution to the Gauss-Bonnet theorem for the incomplete metric $x_{\fun}^2g$ (where $x_{\fun}$ is a bdf for $\pa_{\fun}M$).
\end{remark}

\section{A Polyakov formula for the renormalized determinant}

In this section we will explain how the renormalized integrals described above can be used to define renormalized traces which in turn can be used to
define the determinant of the Laplacian on non-compact surfaces.

\subsection{The renormalized determinant} \label{sec:RDet} $ $\newline
Since the work of Ray and Singer on analytic torsion \cite{Ray-Singer},
zeta regularization has been used to define the determinant of the Laplacian.
Starting with the identity, for $\{\lambda_i \} \subseteq \bbR^+,$
\begin{equation*}
    \pa_s\rest{s=0}\lrpar{ \sum_{i=1}^N \lambda_i^{-s} }
    = - \sum_{i=1}^N \log \lambda_i
    = - \log \prod_{i=1}^N \lambda_i ,
\end{equation*}
Ray and Singer proposed to define the determinant of the Laplacian of $(M,g)$ by first setting
\begin{equation}\label{NaiveZeta}
    \zeta(s)
    = \sum_{\lambda \in \Spec(\Delta)\setminus \{ 0 \}} \lambda^{-s}
\end{equation}
and then formally defining
\begin{equation}\label{FirstDefDet}
    \det \Delta = e^{-\pa_s\rest{s=0}\zeta}.
\end{equation}
To make sense of this formula we note that if $M$ is closed then, by Weyl's law, the sum defining $\zeta(s)$ converges if $\Re(s) > \dim M/2$ and defines a holomorphic function on this half-plane.
It is possible to extend this function meromorphically to the whole plane and then \eqref{FirstDefDet} involves taking the derivative of this meromorphically extended function at the origin, which turns out to be a regular point.

One way to justify the meromorphic extension is to rewrite $\zeta(s)$ using the heat kernel of $\Delta.$
It is easy to check that, for any $\lambda >0,$
\begin{equation*}
    \lambda^{-s} = \frac1{\Gamma(s)} \int_0^{\infty} t^s e^{-t\lambda} \frac{dt}t
\end{equation*}
and hence for $\Re(s) > \dim M/2 ,$
\begin{equation}\label{FirstDefZeta}
    \zeta(s) = \frac1{\Gamma(s)} \int_0^{\infty} t^s \Tr(e^{-t\Delta} - \cP) \frac{dt}t
\end{equation}
where $\cP$ is the projection onto the null space of $\Delta,$ i.e., the constant functions.
As $t \to 0$ the trace of the heat kernel has an asymptotic expansion
\begin{equation*}
    \Tr(e^{-t\Delta}) \sim t^{-\dim M/2} \sum_{k \geq 0} a_k t^k
\end{equation*}
which implies that $\zeta(s)$ has a meromorphic continuation to the complex plane with potential poles at
\begin{equation*}
    s \in \lrbrac{ \frac{\dim M}2, \frac{\dim M}2 - 1, \frac{\dim M}2 - 2, \ldots }.
\end{equation*}
It can then be explicitly checked that zero is a regular point so that the right hand side of \eqref{FirstDefDet} is well-defined.

In terms of the renormalized integrals of the previous section, since
\begin{equation*}
    \frac1{\Gamma(s)} \sim s + \cO(s^2), \Mas s \to 0,
\end{equation*}
we have
\begin{equation*}
    \zeta'(0) = \sideset{^R}{_0^\infty}\int  \Tr(e^{-t\Delta} ) \; \frac{dt}t
\end{equation*}
where $t$ is used to renormalize the integral at $t=0$ and $t^{-1}$ is used to renormalize as $t \to \infty.$
This equality holds on closed manifolds, but notice that the right hand side can serve as a definition of the (logarithm of the) determinant of the Laplacian that does not require that the origin be a regular point for the meromorphically continued zeta function.

On a non-compact surface $(M,g),$ the spectrum of the Laplacian consists of eigenvalues and a continuous spectrum, so one cannot define $\zeta(s)$ by \eqref{NaiveZeta}.
There is also a problem with extending definition \eqref{FirstDefZeta} because the heat kernel of the Laplacian is not of trace class.
However this problem can be overcome by means of a {\em renormalized} trace.

To motivate the renormalized trace recall Lidskii's theorem which says that, if $A$ is an operator that acts via a continuous kernel $\cK_A,$
\begin{equation*}
    Af(\zeta)
    = \int \cK_A(\zeta, \zeta') f(\zeta') \; d\zeta',
\end{equation*}
and $A$ is of trace class, then
\begin{equation*}
    \Tr(A) = \int \cK_A(\zeta, \zeta) \; d\zeta.
\end{equation*}
Since we already know how to renormalize integrals, we define {\bf the renormalized trace} of an operator $A$ to be
\begin{equation*}
    \RTr A = \Rint \cK_A(\zeta, \zeta) \; d\zeta,
\end{equation*}
whenever the right-hand side makes sense.

Fortunately, the heat kernel of a $\Fhc$ metric is well-enough understood to define its renormalized trace.
The heat kernel for a metric with ends asymptotic to hyperbolic cusps is described in \cite{Vaillant} and for a metric with ends asymptotic to hyperbolic funnels in \cite{Albin2}, and it is straightforward to patch together a heat kernel for a general $\Fhc$ metric from these.
In either case it is shown that the distributional kernel of $e^{-t\Delta}\rest{\diag}$ is a smooth function in the interior of $M$ with polyhomogeneous expansions at the boundary of $M$ and also as $t \to 0.$  From the analysis of the Laplacian in \cite{Mazzeo:Hodge} and \cite{Vaillant}, we know that $\Delta$ has closed range and hence $e^{-t\Delta}$ converges exponentially to the projection onto the null space of $\Delta$ as $t \to \infty.$
These properties allow us to make the following definition.

\begin{definition}
Let $M$ be a non-compact surface with a $\Fhc$ metric $g,$ and let $x$ be a total boundary defining function.
The {\bf renormalized trace of the heat kernel} of $g$ is defined to be
\begin{equation*}
    \RTr{e^{-t\Delta}}
    = \Rint \cK_{e^{-t\Delta}}\rest{\diag} \; \dvol_g
\end{equation*}
where the renormalization is carried out using $x.$  From \cite{Albin2} and \cite{Vaillant} (see \cite[Appendix]{Albin-Rochon3}) we know that, as $t \to 0,$
\begin{equation}\label{SmallTimeTr}
    \RTr{ e^{-t\Delta} }
    \sim
    \sum_{k \geq -2} a_{k} t^{k/2}
    + \sum_{k \geq -1} \wt a_k t^{k/2} \log t.
\end{equation}
The {\bf determinant of the Laplacian} is defined to be
\begin{equation*}
    \det \Delta = \exp \lrpar{ - \sideset{^R}{_0^{\infty}} \int  {}^R\!\Tr(e^{-t\Delta} ) \; \frac{dt}t }.
\end{equation*}
\end{definition}

As mentioned above, this definition directly extends the definition from operators with trace-class heat kernel. We will show in \S\ref{sec:Selberg} that it also extends the definition of the determinant via the Selberg zeta function from hyperbolic metrics (cf. \cite{BJP:Selberg}, \cite{TZ1}). Another point in its favor is the main result of this paper: it singles out constant curvature metrics as its critical metrics.

In \eqref{SmallTimeTr}, the logarithmic terms come from the cusp ends.  In fact, as the next lemma shows, many of the coefficients $\wt a_{k}$ automatically vanish. \begin{lemma}
For any $\Fhc$ metric, the terms $\wt a_k$ in the short-time asymptotics of the renormalized trace of its heat kernel \eqref{SmallTimeTr} vanish if $k$ is even.
In particular, $\wt a_{0}=0.$
\end{lemma}

\begin{proof}
Assume without loss of generality that there is only one cusp end.
Recall from \cite[Appendix]{Albin-Rochon3} that the asymptotic expansion \eqref{SmallTimeTr} is derived by analyzing
the restriction of the integral kernel $\cK$ of the heat kernel to the manifold with corners
\begin{equation*}
  \diag_{H}= [M\times \bbR^{+}_{\sqrt{t}}; \pa M\times \{0\}]
\end{equation*}
where it is polyhomogeneous.
This space has three boundary hypersurfaces: $\bhs{\bbf}$ and $\bhs{\tf}$ coming from the `old' hypersurfaces $\{ x= 0 \}$ and $ \{ t = 0 \}$ in $M \times \bbR^+_{\sqrt t}$, respectively, and $\bhs{\cuf}$ the `front face' resulting from blowing-up $\pa M \times \{0\}$.
The logarithmic terms in the expansion come from the corner $\bhs{\cuf} \cap \bhs{\tf}$ and, since we are renormalizing, also from the corner $\bhs{\bbf} \cap \bhs{\cuf}$.

At the latter corner, the log terms can be shown to arise from the short-time expansion of the coefficient of $x^{-1}$ in the expansion of the heat kernel at the face $\bhs{\bbf}$. This  can be identified following Vaillant \cite[Chapter 4]{Vaillant} with the heat kernel of a model operator (the `horizontal family') which in this case  is multiplication by a constant, namely $-\frac r8$ where $r$ is the asymptotic value of the scalar curvature at the cusp. Thus this corner contributes $\frac{ e^{ \frac r8 t} }{\sqrt t} \log t$ to the short-time expansion of the renormalized heat trace, and so does not contribute to $\wt a_k$ for $k$ even.

To handle the former corner, $\bhs{\cuf} \cap \bhs{\tf}$, recall the well-known fact that for closed manifolds of dimension $n$, the short-time asymptotics of the heat trace are of the form $t^{-n/2}$ times an expansion in $t$ (as opposed to $\sqrt t$). The same `even-ness' is true in the expansion of the heat kernel of a $\Fhc$ metric at the boundary face $\bhs{\tf}$. It follows (see \cite[Appendix]{Albin-Rochon3}) that the asymptotics at the corner $\bhs{\cuf} \cap \bhs{\tf}$ do not contribute to $\wt a_k$ for $k$ even.
\end{proof}

We can also define the $\zeta$ function of the Laplacian of $g,$ for $\Re s \gg 0,$ by
\begin{equation}\label{RenZetaFun}
    \zeta(s) = \frac1{\Gamma(s)} \int_0^{\infty} t^s \; \RTr{e^{-t\Delta} - \cP} \frac{dt}t
\end{equation}
where $\cP$ is the projection onto the $L^{2}$ null space of $\Delta$ (i.e., the constant functions if $g$ has finite area and the zero function otherwise).
Using \eqref{SmallTimeTr} this function has a meromorphic continuation to the complex plane, still denoted $\zeta.$
Since $\wt a_0$ vanishes, zero is a regular point of $\zeta(s)$ and the definition above coincides with
\begin{equation*}
    \det\Delta = e^{-\zeta'(0)}.
\end{equation*}

\begin{remark}
The renormalized trace was first defined by Melrose in his proof of the Atiyah-Patodi-Singer index theorem \cite{APSBook} for asymptotically cylindrical metrics.
Melrose also pointed out that this definition allows one to extend the definition of the zeta function and the determinant of the Laplacian as described above.
See \cite{Hassell} for an analysis of this extension to asymptotically cylindrical metrics including a proof of the Cheeger-M\"uller theorem for the corresponding analytic torsion.
\end{remark}

\subsection{Relation with the relative determinant} $ $\newline
A circle of ideas from scattering theory was introduced into this context by
Werner M\"uller \cite{Mueller:RD} as a way to overcome the problem of the heat
kernel of the Laplacian not being of trace class. Instead of extending the trace
functional to operators that are not of trace class a model operator is
introduced so that the difference of the two operators is of trace class.
One can consider an operator that is naturally related to the
surface. In the case of finite area, i.e., only cusps ends, the
natural model operator is the following: Decompose $M$ as a compact
part, $M_{0},$ and the cusps ends, $Z_{hc}.$ The cusps ends are each
considered with one boundary that joins the end to the compact
part $M_{0}.$
The model operator is the direct sum of the
self-adjoint extension of the hyperbolic Laplacian on the cusp ends
with respect to Dirichlet boundary conditions at the boundaries, and
the operator zero on $M_{0}.$ We denote this operator as
$\Delta_{0}.$ The definition of the relative determinant depends on the following facts:

\begin{enumerate}
\item $e^{-t \Delta_{g}} - e^{-t \Delta_{0}}$ is  an
operator of trace class for all $t>0.$
\item As $t\to 0^{+},$ there is an asymptotic expansion of the relative heat trace of the form given in
\eqref{SmallTimeTr} with $\wt a_0 =0$.
\item Since there is a spectral gap at zero we have that:
$$\Tr(e^{-t \Delta_{g}}-e^{-t \Delta_{0}}) = 1 + O(e^{-ct}),$$
as $t\to \infty,$ where $1=\dim \ker \Delta_{g} -\dim \ker
\Delta_0.$
\end{enumerate}
Then we can define the relative zeta function as:
$$\zeta(s;\Delta_g,\Delta_0)= {\frac{1}{\Gamma(s)}} \int_{0}^{\infty}
t^s(\Tr(e^{-t \Delta_g}-e^{-t \Delta_{0}})-1)  \frac{dt}t.$$ In the same
way as for the regularized determinant, the asymptotic expansions
guarantee that the relative zeta function has a meromorphic
continuation to the complex plane that is regular at zero, and we
can define: $$\det(\Delta_g, \Delta_0) =
e^{-\zeta'(0;\Delta_g,\Delta_{0})}.$$
The relative determinant has been considered on non-compact surfaces where the metric is exactly hyperbolic outside a compact set by M\"uller \cite{Muller:Inv}, \cite{Mueller:RD} and by Borthwick, Judge, and Perry \cite{BJP:Duke}.
In a more recent work \cite{Aldana}, the second author has extended the definition of the relative determinant to surfaces whose ends are merely asymptotic to hyperbolic cusps.

In either case, working with renormalized traces allows us to write
\begin{gather*}
    \Tr(e^{-t \Delta_{g}}-e^{-t \Delta_{0}})
    = \RTr{e^{-t\Delta_g}} - \RTr{e^{-t\Delta_0}}, \\
    \zeta(s;\Delta_g,\Delta_0) = \zeta(s;\Delta_g) - \zeta(s; \Delta_0), \Mand
    \det(\Delta_g, \Delta_0) = \frac{\det(\Delta_g)}{\det(\Delta_0)}.
\end{gather*}
Hence, when the relative determinant and the renormalized determinant are both defined, they are essentially equivalent.

\subsection{Relation with Selberg zeta function} \label{sec:Selberg} $ $ \newline
Both in order to compute the critical value of the determinant of the Laplacian, and to relate to other approaches
to the determinant of non-compact surfaces, it is important to understand the relation between the determinant as defined above
and the Selberg zeta function.

For compact hyperbolic Riemann surfaces, an interesting relationship between the determinant
of the Laplacian and the Selberg Zeta function was discovered by D'Hoker and Phong \cite{DhokerPhong} and further refined by Sarnak \cite{Sarnak}.
If $M$ is a non-compact surface with a $\Fhc$-metric $g$ of {\em constant curvature}, then this relationship still holds.
We can write $M$ as $\bbH/\Gamma$ and define the Selberg zeta function $Z(s),$ for $\Re s >1,$ by the absolutely convergent product
\begin{equation*}
    Z(s) = \prod_{ \{ \gamma \} } \prod_{k=0}^{\infty} (1 - e^{-(s+k) \ell(\gamma) } )
\end{equation*}
where the outer product goes over conjugacy classes of primitive hyperbolic elements of $\Gamma,$ and $\ell(\gamma)$ is the length of the corresponding closed geodesic.
The function $Z(s)$ admits a meromorphic continuation to the whole complex plane.

Borthwick, Judge, and Perry showed \cite[(5.2)]{BJP:Selberg} that if
\begin{equation*}
    R_g(s) = (\Delta + s(s-1))^{-1}
\end{equation*}
and $\cL(s)$ is any function satisfying
\begin{equation}\label{BJPEq}
    \lrpar{ \frac1{2s-1} \frac{\pa}{\pa s}}^2 \log \cL(s) = - \RTr{ R_g(s)^2 }
\end{equation}
then there are constants $E$ and $F$ such that
\begin{equation}\label{BJPCon}
    \cL(s) = Z(s) e^{E + Fs(1-s)} \lrpar{ \frac{\Gamma(s)}{(2\pi)^s \Gamma_2^2(s) } }^{\chi(M)}
\lrpar{ 2^{s} \sqrt{\pi(s-\tfrac12)} \Gamma(s-\tfrac12)  }^{-n_C}
\end{equation}
(here $\Gamma_2(s)$ is Barnes' double Gamma function).
We will show that, with the determinant defined above, $\det( \Delta + s(s-1) )$ satisfies \eqref{BJPEq} and hence \eqref{BJPCon}. As in \cite{Sarnak}, by examining the asymptotics as $s\to\infty$ we will determine the values of $E$ and $F$ in this case.

\begin{theorem}\label{thm:RelSelberg}
Let $M$ be a non-compact surface and $g$ a $\Fhc$-metric on $M$ of constant curvature.
The zeta regularized determinant $\det( \Delta + s(s-1) )$ satisfies \eqref{BJPEq} and hence \eqref{BJPCon}.
The constants $E$ and $F$ in this case are equal to
\begin{equation*}
    E= \chi(M) \lrpar{ \frac12\log 2\pi - 2\zeta_R'(-1) + \frac14} , \quad
    F=-\chi(M)
\end{equation*}
where $\zeta_R$ is the Riemann zeta function.
It follows that
\begin{equation*}
    \det (\Delta) = \begin{cases}
        C_{\Fhc} \; Z'(1) & \text{ if Area}(g) < \infty \\
        C_{\Fhc} \; Z(1) & \text{ otherwise}.
        \end{cases}
\end{equation*}
with
\begin{equation*}
    C_{\Fhc}= e^E (2\pi)^{-\chi(M)}(\sqrt 2\pi)^{-n_C}.
\end{equation*}
\end{theorem}

\begin{remark}
The reason a derivative is taken in the case of a finite area surface is that to compute $\det(\Delta)$ one needs to exclude the zero eigenvalue of $\Delta.$
\end{remark}

\begin{proof}
In \cite[(7.9)]{Albin-Rochon2} it is shown that \eqref{SmallTimeTr} implies
\begin{equation}\begin{split}
    -\log\det\lrpar{\Delta + w}
    &=\int_0^\infty  \lrpar{{}^R\!\Tr\lrpar{e^{-t\Delta} }-f_0\lrpar t} e^{-tw} \;\frac{dt}t\\
    &-a_0\log w
    -2\sqrt\pi a_{-1}\sqrt w
    +a_{-2}w\lrpar{-1+\log w} \\
    &+\wt a_{-1}\sqrt{w}\left( \Gamma_{\log}(-\frac12)- \log w \Gamma(-\frac12)\right)
\end{split}
\label{ExplicitMero}\end{equation}
where the $a_{k}$ and $\wt a_{k}$ are the coefficients in the short time asymptotic of the trace of the heat kernel \eqref{SmallTimeTr} and
\begin{equation}\label{gamma_log.1}
\begin{gathered}
    f_0(t) = a_{-2}t^{-1} + \wt a_{-1} t^{-1/2} \log t + a_{-1} t^{-1/2} + a_0, \\
    \Mand
    \Gamma_{\log}(z):= \int_{0}^{\infty} t^{z} e^{-t}\log t \frac{dt}{t}.
\end{gathered}
\end{equation}
If we set $w=s(s-1)$ so that $\frac{\pa}{\pa w}= \frac{1}{2s-1}\frac{\pa}{\pa s}$,
it follows that (\cite[(7.21)]{Albin-Rochon2})
\begin{multline*}
    \lrpar{ \frac1{2s-1} \frac{\pa}{\pa s} } \log\det (\Delta + s(s-1)) \\
    =\int_0^\infty  \lrpar{{}^R\!\Tr\lrpar{e^{-t\Delta}}-a_{-2}t^{-1}} e^{-ts\lrpar{s-1}} \;dt
    -a_{-2}\log \lrpar{ s(s-1)}
\end{multline*}
and hence
\begin{multline*}
    \lrpar{ \frac1{2s-1} \frac{\pa}{\pa s} }^2 \log\det (\Delta + s(s-1))
    =-\int_0^\infty  t \left(\RTr{e^{-t\Delta}}\right) e^{-ts(s-1)} \; dt \\
    = -\RTr{ \int_0^\infty t e^{-t(\Delta + s(s-1))} \; dt }
    = - \RTr{ (\Delta +s(s-1) )^{-2} }.
\end{multline*}
which establishes \eqref{BJPEq} and \eqref{BJPCon}.

To determine $E$ and $F,$ we can then proceed exactly as in the
proof of \cite[Theorem 2]{Albin-Rochon2}.  Finally, because zero
is in the spectrum of the Laplacian on a hyperbolic surface
precisely when its area is finite,
\begin{equation*}
    \det (\Delta) = \begin{cases}
    \lim_{s \to 1}\frac{\det( \Delta + s(s-1) )}{s(s-1)}  & \text{ if Area}(g) < \infty \\
    \det( \Delta + s(s-1) )\rest{s=1} & \text{ otherwise }
    \end{cases}
\end{equation*}

\end{proof}

In fact, the proof of \cite[Theorem 2]{Albin-Rochon2} also give us
the values of the first few coefficients in the short time
asymptotic \eqref{SmallTimeTr} of the trace of the heat kernel
\begin{equation}
\begin{gathered}
    a_{-2} = \frac{{}^{R}\!\!\Area(M)}{4\pi}, \quad
    \wt a_{-1} = \frac{n_C}{4\sqrt \pi}, \quad
    a_0 = \frac{\chi(M)}6, \\
    a_{-1}
    =  \frac{n_{C}}{2\sqrt\pi}\lrpar{
    \frac{\Gamma_{\log}(-\tfrac12)}{4\sqrt \pi} +1-\log 2 }.
\end{gathered}
\label{coef.1}\end{equation}

When the scalar curvature is not constant but the metric has totally
geodesic ends, we get the same coefficients as the next lemma shows.

\begin{lemma}\label{lem:HeatKer}
Let $g$ be a $\Fhc$ metric on $M.$

\begin{itemize}
\item [1)]
If the ends of $g$ are totally geodesic with scalar curvature asymptotically equal to $-2$ in each end, then in the expansion \eqref{SmallTimeTr} the first coefficients are given by \eqref{coef.1}.

\item [2)]
If $\varphi = \cO(x^2_{\fun})$ and $\varphi=\cO(x_{\hc}),$ then $\int_0^t e^{-(t-s)\Delta}\varphi\Delta e^{-s\Delta} \; ds$ is of trace class.
\end{itemize}
\end{lemma}

\begin{remark}
Although we only compute the coefficients $a_{-1}$ and $\wt a_{-1}$ when the
scalar curvature is equal to $-2,$ they can be obtained in the general case by
an appropriate rescaling argument.
\end{remark}

\begin{remark}
If there are only cusps,  see also \cite{Aldana}  for a different proof of ($2$).
\end{remark}

\begin{proof}
For the first part of the lemma, notice that the `interior contributions' are
given by the same integration of local terms, but with integrals replaced with renormalized integrals whenever the volume is infinite.  Thus,
we need to check the contributions coming from the cusp and funnel ends are the same for a metric with totally geodesic ends and a metric which is hyperbolic
in a collar neighborhood of each end.

First let us focus on what happens near a cusp.  Without loss of generality, we can assume that we have only cusp ends, in fact only one cusp.  By assumption, we know that
\eqref{ModelCusp} is satisfied with $\left. \varphi\right|_{\pa \bar M}=0.$
For these metrics the heat kernel is an element of Vaillant's heat calculus \cite[\S 4]{Vaillant},
\begin{equation*}
    e^{-t\Delta_g} \in \Psi_H^{2,2,0}(M).
\end{equation*}
The superindices in $\Psi^{2,2,0}_H(M)$ refer to the behavior of $e^{-t\Delta}$ as $t\to 0$ in the interior of $M,$ as $t\to 0$ and $x\to 0,$ and as $x \to 0$ for $t>0.$  Let $g_{\hc}$ be a metric which is hyperbolic in a neighborhood of the cusp end.
We will see, using Duhamel's formula and Vaillant's composition formula, that the difference between $e^{-t\Delta_g}$ and the heat kernel of $g_{\hc}$ vanishes at the cusp as $t \to 0.$
Suppose that in a neighborhood $E$ of $\pa_{\hc}M$ we can write the metric as
\begin{equation*}
    g = e^{\varphi} \left( \frac{dx^2}{x^2}+ x^2d\theta^2  \right) = e^{\varphi} g_{\hc}
\end{equation*}
for some smooth $\varphi = \cO(x^k)$ with $k>0.$
By introducing the interpolating family
\begin{equation*}
    g_\tau = e^{\tau \varphi} g_{\hc},
\end{equation*}
the heat kernels of $g$ and $g_{\hc}$ in this neighborhood are related by
\begin{equation*}
\begin{gathered}
    e^{- t\Delta_g}
    = e^{-t\Delta_{\hc}}
    + \int_0^1 \pa_\tau e^{-t\Delta_\tau} \; d\tau \\
    = e^{-t\Delta_{\hc}}
    - \int_0^1
    \int_0^t e^{-(t-s)\Delta_\tau} \varphi \Delta_\tau e^{-s\Delta_\tau} \; ds \; d\tau.
\end{gathered}
\end{equation*}
From the proof of \cite[Theorem 4.9]{Vaillant}, we know that $\Delta_\tau e^{-s\Delta_\tau}$ is an element of $\Psi_{H}^{0,0,0}(M)$ in Vaillant's heat calculus.
For $\varphi$ in $\cO(x^k),$ $\varphi \Delta_\tau e^{-s\Delta_\tau}$ is an element of $\Psi_{H}^{0,k,k}(M)$ and hence using Vaillant's composition formula \cite[Theorem 4.6]{Vaillant}
\begin{equation*}
    \int_0^t e^{-(t-s)\Delta_\tau} \varphi \Delta_\tau e^{-s\Delta_\tau} \; ds \in \Psi_H^{2,2+k,\cH}(M)
\end{equation*}
for some index set $\cH$ bounded below by $k.$
This implies that the $a_{j}$ and $\wt a_j$  in the short-time asymptotics \eqref{SmallTimeTr} are given by the same formula as those in the short-time asymptotics of $e^{-t\Delta_{\hc}}$ for $j< (k-1)/2.$
Elements in $\Psi_H^{2,2+k,\cH}(M)$ are trace-class at positive time for $k >0,$ so the discussion above also establishes ($2$) along the cusp ends. A similar argument using \cite{Albin2} instead of \cite{Vaillant} (and the fact that $\varphi = \cO(x_{\fun}^2)$) establishes the lemma along the funnel ends.

\end{proof}

\subsection{Polyakov formula} \label{sec:Polyakov} $ $\newline
In this section we extend Polyakov's formula for the change in the determinant of the Laplacian upon a conformal change of metric. For the relative determinant, this formula is due to Borthwick, Judge, and Perry \cite{BJP:Duke} (for metrics with no cusps and exactly hyperbolic outside a compact set) and to the second author \cite{Aldana} (for metrics with no funnels).

We will assume that the metrics involved are $\Fhc$ metrics with totally geodesic ends.
As explained in the previous section, the latter assumption simplifies the behavior of the heat kernel and includes the metrics of constant scalar curvature.
Thus we will analyze the behavior of the determinant of the Laplacian for a family of metrics $g(\tau) = e^{\omega(\tau)}g_0$ where $g_0$ is a smooth $\Fhc$ metric and
\begin{equation}\label{OmegaCond}
    \omega(\tau) = \wt \omega(\tau)+ \sum_{i=1}^{n_{\fun}+n_{\hc}} \omega_{i}(\tau)\chi(x_{i}), \Mwith \wt\omega \in x^{2}\CI(\bar M\times [0,T]_{\tau}),
\end{equation}
where $\omega_{i}\in\CI(\bbR)$ and $\chi\in\CI_{c}([0,+\infty)_{u})$ is a cut-off
function equal to 1 for $u<\frac{\epsilon}{2}$ and to $0$ for $u>\frac{3\epsilon}{4}$.

\begin{theorem}[Polyakov formula] \label{thm:Polyakov}
Let $(M,g_0)$ be a non-compact surface with a smooth $\Fhc$ metric.
Let $\omega(\tau) \in \CI(\bar M)$ be a smooth family of functions satisfying \eqref{OmegaCond} and $g(\tau) = e^{\omega(\tau)}g_0$.   Then the determinant of the Laplacian satisfies
\begin{equation}\label{Polyakov}
    \pa_{\tau} \log \det \Delta_\tau
    = - \frac1{24\pi} \overset{R\;\;\;}{\int_{M}} \omega'(\tau) R_\tau \dvol_\tau
\end{equation}
if the area of $M$ is infinite, and
\begin{equation}\label{PolyakovFinite}
    \pa_{\tau} \log \det \Delta_\tau
    =
    -\frac1{24\pi} \int_M \omega'(\tau) R_\tau \dvol_{\tau}
    + \pa_\tau \log \Area_\tau(M)
\end{equation}
if the area of $M$ is finite.
\end{theorem}

\begin{proof}
Let $\Delta_\tau$ denote the Laplacian of $g(\tau),$ so that $\Delta_{\tau} = e^{-\omega(\tau)}\Delta_0.$
Then
\begin{equation*}
\begin{split}
    \pa_{\tau} {}^R\!\Tr (e^{-t\Delta_{\tau}})
    &= \FP_{z =0} \pa_{\tau} \Tr (x^z e^{-t\Delta_{\tau}}) \\
    &= - \FP_{z=0} \Tr \lrpar{ x^z \int_0^t e^{-(t-s)\Delta_\tau} (-\omega'(\tau)) \Delta_\tau e^{-s\Delta_{\tau}} \; ds } \\
    &= \sum_{i=1}^{n_{\fun}+n_{\hc}} \omega_i'(\tau) \FP_{z=0}
     \Tr \lrpar{ x^z t \chi(x_{i})\Delta_\tau e^{-t\Delta_\tau} }   +  \Tr \lrpar{ \wt\omega'(\tau) t \Delta_\tau e^{-t\Delta_\tau} } \\
     & + \sum_{i=1}^{n_{\fun}+n_{\hc}} \omega_i'(\tau)
     \int_{0}^{t} \FP_{z=0} \Tr\lrpar{x^{z} [e^{-(t-s)\Delta_\tau}, \chi(x_{i}) \Delta_\tau e^{-s\Delta_{\tau}}] } \;ds   \\
     &=  \RTr{ \omega'(\tau) t \Delta_\tau e^{-t\Delta_\tau} }
    = -t \pa_t \RTr{ \omega'(\tau) e^{-t\Delta_\tau} }
\end{split}
\end{equation*}
where in the third equality we have used that $\int_0^t e^{-(t-s)\Delta_\tau} \wt \omega'(\tau) \Delta_\tau e^{-s\Delta_{\tau}}$ is of trace class and that
\[
    \RTr{ [e^{-(t-s)\Delta_\tau}, \chi(x_{i}) \Delta_\tau e^{-s\Delta_{\tau}}] }=0.
\]
Indeed, since this term is the regularized trace of a commutator, its value
depends on the asymptotic expansion of the various operators at $Y_{i}$.
Since $\chi(x_{i})\equiv 1$ near $Y_{i}$, this means
\[
\RTr{ [e^{-(t-s)\Delta_\tau}, \chi(x_{i}) \Delta_\tau e^{-s\Delta_{\tau}}] }=
\RTr{ [e^{-(t-s)\Delta_\tau},\Delta_\tau e^{-s\Delta_{\tau}}]}=0,
\]
the latter commutator vanishing identically.

Next let $\cP = \cP(\tau)$ denote the $L^2$ projection onto constants in the case of finite volume and zero otherwise.
We will use the obvious facts that $\pa_\tau \Tr(\cP) =0$ and $\pa_t \Tr (\omega'(\tau)\cP) =0.$

We can proceed as follows
\begin{equation*}\label{FirstZetaStep'}
\begin{split}
    \ddt{\tau} \dzero{s} \zeta_{e^{\omega(\tau)}g}(s)
    &= \dzero{s} \lrpar{ \frac1{\Gamma(s)} \int_0^{\infty} t^s \ddt{\tau} {}^{R}\!\Tr ( e^{-t\Delta_{\tau}} - \cP ) \frac{dt}t } \\
    &= \dzero{s} \left( -\frac{1}{\Gamma(s)}
    \int_0^\infty t^s \pa_t {}^R\!\Tr( \omega'(\tau) (e^{-t\Delta_{\tau}} - \cP) ) \; dt
    \right) \\
    &= \dzero{s} \left( \frac{s}{\Gamma(s)}
    \int_0^\infty t^{s-1} \; {}^R\!\Tr( \omega'(\tau) (e^{-t\Delta_{\tau}} - \cP) ) \; dt \right)
\end{split}
\end{equation*}
Since $\frac1{\Gamma(s)} \sim s + \cO(s^2)$ as $s\to0,$ this is equal to
\begin{equation}\label{InteriorTerm1}
    \Res_{s=0}
    \int_0^\infty t^s \; {}^{R}\!\Tr( \omega'(\tau) (e^{-t\Delta_{\tau}} - \cP) ) \; \frac{dt}t
    =  \FP_{t=0} {}^{R}\!\Tr( \omega'(\tau) (e^{-t\Delta_{\tau}} - \cP) )
\end{equation}

In the interior, we know that
$\omega'(\tau)e^{-t\Delta}$ has a short-time expansion of the form
\begin{equation*}
\omega'(\tau) (e^{-t\Delta_{\tau}} - \cP) \sim \frac{\alpha_{-1}}{t}
+\alpha_{0} + o(1)
\end{equation*}
where the coefficients $\alpha_k$ are precisely those functions that occur on a closed surface, multiplied by $\omega'(\tau),$ e.g. $\alpha_0 = \frac R{24\pi} \omega'(\tau).$
Because $\omega'(\tau)-\omega_{i}'(\tau)= \mathcal{O}(x^{2}_{i}),$ we then see from Lemma~\ref{lem:HeatKer} as well as the construction of the heat kernel in \cite{Vaillant} for cusps (see \cite[Appendix]{Albin-Rochon3}) and \cite{Albin2} for funnels that
\begin{equation*}
\FP_{t =0}\Tr ( \omega'(\tau)e^{-t\Delta_{\tau}} )
    = \frac1{24\pi} \int_M \omega'(\tau) R_{\tau} \dvol_{\tau}.
\end{equation*}
Thus, we conclude that
\begin{gather*}
     \FP_{t =0} {}^{R}\!\Tr ( \omega'(\tau) (e^{-t\Delta_{\tau}} - \cP) ) \\
    = \frac1{24\pi} \int_M \omega'(\tau) R_{\tau} \dvol_{\tau}
    - \frac {1}{\Area_\tau(M)} \int_M \omega'(\tau) \dvol_{\tau}
\end{gather*}
where the final term is replaced by zero if the volume of $M$ is infinite, and, if not, can be rewritten
\begin{equation*}
    \frac {1}{\Area_\tau (M)} \int_M \omega'(\tau) \dvol_{\tau}
    = \frac1{\Area_\tau (M)} \int_M \pa_\tau (e^{\omega(\tau)}) \dvol_0
    = \pa_\tau \log \Area_\tau (M).
\end{equation*}
Finally, since
\begin{equation*}
    \ddt{\tau} \dzero{s} \zeta_{e^{\omega}g}(s) = -\ddt{\tau} \log\det\Delta_\tau
\end{equation*}
this finishes the proof.
\end{proof}

To study the determinant of the Laplacian, it is natural to impose the following
normalization condition
\begin{equation}
    {}^{R}\!\!\Area(g)= -2\pi \chi(M),
\label{norm.1}\end{equation}
as otherwise one can increase the determinant `artificially' by scaling the metric (e.g., consider $g(\tau) = e^{\tau C}g$ for a constant $C$ in the Polyakov formula).  Notice that the case where $\chi(M)=0$ is special, since \eqref{norm.1} then implies there must be at least one funnel end, but it does not rule out this sort of scaling.  On the other hand, in this case, these scalings  are no longer problematic, since they leave the determinant invariant as can be seen from the Polyakov formula.   This leads naturally to the following definition.

\begin{definition}
A $\Fhc$ metric $g$ satisfying \eqref{norm.1} is said to be critical for the determinant of the Laplacian
if for any $\psi\in \CI_{c}(M)$ with $\int_{M} \psi dg=0,$ we have that
\[
        \int_{M} \psi R_{g} dg=0.
\]
Thus, when $\chi(M)\ne 0$, a critical $\Fhc$ metric is precisely one with constant scalar curvature $-2$.  When instead $\chi(M)=0$, a critical $\Fhc$ metric is one with constant negative scalar curvature.
\label{critical.1}\end{definition}

Now suppose that $\omega(\tau) = \tau \omega$ for a fixed smooth function $\omega$ satisfying
\begin{equation}\label{VanishingAss}
    \omega = \omega_0 + \wt\omega, \Mwith \wt\omega =\cO(x^2) \quad \mbox{and} \quad
    \omega_{0} \; \mbox{a constant},
\end{equation}
and let $g_1 = e^{\omega}g_0.$
For metrics of finite area we can integrate equation \eqref{PolyakovFinite} and use \eqref{conf.2a} to find
\begin{equation}\label{IntegratedPolyakov}
\begin{gathered}
    \log \det \Delta_{g_1} - \log \det \Delta_{g_0}
    = \int_0^1 \pa_\tau \log\det \Delta_{e^{\tau\omega}g_0} \; d\tau \\
    =\int_0^1 \lrspar{
    -\frac1{24\pi} \int_M \omega' R_\tau \dvol_{\tau}
    + \pa_\tau \log \Area_\tau(M) } \; d\tau \\
    =  \log \Area_1(M) - \log \Area_0(M)
    -\frac1{24\pi} \int_M \lrpar{ \omega' R_0 + \tfrac12|\nabla_0\omega|^2} \dvol_0
\end{gathered}
\end{equation}
where \eqref{VanishingAss} guarantees both that there is no boundary term from applying Green's theorem and that $|\nabla\omega|^2$ is integrable.
In the same way, for infinite area metrics, integrating \eqref{Polyakov} yields
\begin{equation}
\begin{gathered}
    F(\omega)= \log \det \Delta_{g_1} - \log \det \Delta_{g_0}
    =
    -\frac1{24\pi} \Rint \lrpar{ \omega R_0 + \frac12 |\nabla_0\omega|^2} \dvol_0 \\
    = -\frac{\chi(M)}{6}\omega_0
    -\frac1{24\pi} \int \lrpar{ \wt\omega' R_0 + \tfrac12|\nabla_0\wt\omega|^2} \dvol_0.
\end{gathered}
\label{convex.1}\end{equation}

\section{Ricci Flow on surfaces with funnel, cusp metrics} \label{sec:Ricci}

Hamilton \cite{Hamilton}  (see also \cite{Cao}) studied the Ricci flow
on closed surfaces and showed that, if the Euler characteristic is
negative, then a solution to the normalized Ricci flow exists for
all time and converges exponentially to a hyperbolic metric in the
conformal class of the original metric. Hamilton's result and
approach were extended to non-compact surfaces with asymptotically
hyperbolic cusp ends by Ji, Mazzeo, and \v{S}e\v{s}um
\cite{Ji-Mazzeo-Sesum}. In this section, we further extend their
result to non-compact surfaces whose ends are asymptotic to
funnels or hyperbolic cusps.

In this section we will assume less regularity on the metrics we work with than in the previous sections.
Let $M$ be a surface with boundary and choose a background $\Fhc$ metric $\df g$ that is exactly hyperbolic in a neighborhood of the cusp ends.
For $\alpha \in (0,1)$ and a continuous function $v$ define
\begin{equation*}
    \norm{v}_{0,\alpha} = \sup_{\zeta \in M} |v(\zeta)|
   + \sup \lrbrac{ \frac{|v(\zeta) - v(\zeta')|^{\alpha}}{d(\zeta, \zeta')} : d(\zeta, \zeta') <1 },
\end{equation*}
where the distance between two points is measured with respect to $\df g.$
Let $\cC^{0,\alpha}_{\Fhc}(M)$ denote the space of functions for which $\norm{v}_{0,\alpha} < \infty,$ where along the cusps
the collapse of the injectivity radius is dealt with as in \cite{Ji-Mazzeo-Sesum} by passing to a covering space.
For $k \in \bbN,$ we say that $v\in \cC^{0,\alpha}_{\Fhc}(M)$ is an element of the H\"older space $\cC^{k,\alpha}_{\Fhc}(M)$ if, whenever $V_1, \ldots, V_k$ are vector fields of bounded pointwise length with respect to $\df g,$ we have
\begin{equation*}
    V_1 \cdots V_j v \in \cC^{0,\alpha}_{\Fhc}(M) \quad \mbox{for} \; j\le k.
\end{equation*}
Notice that, along a cusp, vector fields of bounded pointwise length are linear combinations of the vector fields
\begin{equation*}
    x\pa_x, \quad \frac1x\pa_{\theta}.
\end{equation*}
Similarly, we say we say that $v\in \cC^{0}_{\Fhc}(M)$ is an element of the  space $\cC^{k}_{\Fhc}(M)$ if, whenever $V_1, \ldots, V_k$ are vector fields of bounded pointwise length with respect to $\df g,$ we have
\begin{equation*}
    V_1 \cdots V_j v \in \cC^{0}_{\Fhc}(M) \quad \mbox{for} \; j\le k.
\end{equation*}
Both $\cC^{k,\alpha}_{\Fhc}(M)$ and $\cC^{k}_{\Fhc}(M)$ are Banach spaces in a natural way.   Although we will use the H\"older space $\cC^{k,\alpha}_{\Fhc}(M)$ later on (starting in Proposition~\ref{le.1}), we will only need for the moment the simpler space $\cC^{k}_{\Fhc}(M)$.

In this section, we will work with metrics $g_0$ as in Definition~\ref{ModelFun}.  However, instead of assuming the function $\varphi$ is smooth \textbf{up to the boundary}, we will assume that $\varphi\in \mathcal{C}^{k+2}_{\Fhc}(M)\cap \CI(M)$ for some
$k\ge 1.$  We will also assume that for each $i\in\{1,\ldots, n_{\fun}+n_{\hc}\},$ there is a constant $\phi_{i}\in\bbR$ and $\delta>0$ such that
\begin{equation}
\varphi-\phi_{i}\in x_{i}^{\delta}\mathcal{C}^{k+2}_{\Fhc}(M).
\label{se.3b}\end{equation}
It will sometimes be convenient to take different decay conditions for cusp and
funnel ends.  In that case, we will use the notation $\delta=\delta_{\fun}$ for
$i\in\{1,\ldots,n_{\fun}\}$ and $\delta= \delta_{\hc}$ for $i\in\{n_{\fun}+1,\ldots, n_{\fun}+n_{\hc}\}.$

We point out that under these assumptions, the scalar curvature of $g_0$ satisfies
\begin{equation}\label{se.7}
             R_{g_{0}}-r_{i} \in x^{\delta}_{i} \cC^{k}_{\Fhc}(M)
\end{equation}
where $r_{i}= -2e^{-\phi_{i}}$ is the asymptotic value of the curvature
at the end $Y_{i}.$

\subsection{Ricci flow and renormalized area} $ $\newline
Let $M$ be a surface, $g$ a metric on $M$ and $R$ the scalar curvature of $g.$
On a surface the curvature is determined by the scalar curvature, in particular the Ricci curvature of $g$ is $\frac12 R g,$ and the normalized Ricci flow equation is
\begin{equation}\label{RicciFlowEq}
    \begin{cases}
        \pa_t g(t) = (\sC-R_t) g(t) \\ g(0) = g_0
    \end{cases}
\end{equation}
where $\sC$ is a constant.
This flow preserves the conformal class of $g_0,$ so we can write
\begin{equation*}
    g(t) = e^{\omega(t)}g_0
\end{equation*}
for some smooth function $\omega$ which satisfies
\begin{equation}\label{omegaEq}
    \omega'(t) = \sC - R_t.
\end{equation}

It is useful to recall that under a conformal change of metric we have
\begin{equation}
    \Delta_{g(t)} = e^{-\omega(t)}\Delta_0,
    \quad
    \dvol_t = e^{\omega(t)} \dvol_0,
    \quad
    R_{g(t)} = e^{-\omega(t)}(R_{g_{0}} + \Delta_{g_{0}} \omega(t) ),
\label{conf.2a}\end{equation}
where $\Delta_{g(t)}$ is the positive definite Laplacian,
as then from \eqref{omegaEq} it is easy to derive equations for the evolution of these quantities,
\begin{gather}
    \pa_t \Delta_{g(t)} = (R_{g(t)} - \sC) \Delta_t, \quad
    \pa_t \dvol_t = (\sC - R_{g(t)}) \dvol_t, \notag \\
    \pa_t R_t
    = -\Delta_{g(t)} R_{g(t)} + R_{g(t)}(R_{g(t)} - \sC).
    \label{EvolR}
\end{gather}

In particular, the normalized Ricci flow can be written as a scalar equation
\begin{equation}
   \frac{\pa \omega}{\pa t} = - e^{-\omega}(\Delta_{g_{0}}\omega+ R_{g_{0}}) +\sC,
\quad \omega(0)\equiv 0.
\label{se.1}\end{equation}

As on a compact surface, one natural choice for the constant $\sC$ is to take the (renormalized) average curvature.
If the funnel ends are totally geodesic (i.e., $\delta_{\fun}> 1$), so that the renormalized Gauss-Bonnet theorem holds, then
the flow with this choice of $\sC$ will preserve the renormalized area.

\begin{lemma}\label{lem:Vol}
Suppose that $M$ is a non-compact surface and $g(t)$ is a smooth family of
metrics satisfying \eqref{ModelFun} and \eqref{ModelCusp} for some $\varphi\in\mathcal{C}^{k+2}_{\Fhc}(M)$ satisfying \eqref{se.3b} with $\delta_{\fun}>1.$
If we assume that
${}^R\!\!\Area_0(M) \neq 0$ and that
$\pa_t g(t) = (\sC - R_t) g(t)$ with
\begin{equation*}
    \sC =
    \bar R = \frac{ 4\pi \chi(M) }{ {}^R\!\!\Area_0(M) },
\end{equation*}
then, for all $t,$ we have
\begin{equation*}
    {}^R\!\!\Area_t(M) = {}^R\!\!\Area_0(M).
\end{equation*}

If instead we assume that ${}^R\!\!\Area_0(M)= 0$ then, for any $\sC \in \bbR,$
a smooth solution $g(t)$ to $\pa_t g(t) = (\sC - R_t) g(t)$ satisfies ${}^R\!\!\Area_t(M)= 0$ for all $t$
if $\chi(M)=0$, or else $\left. \frac{\pa}{\pa t} {}^{R}\!\!\Area_t(M)\right|_{t=0}\ne 0$ if $\chi(M)\ne 0$.
\end{lemma}

\begin{proof}
Working with an arbitrary value of $\sC,$ we find
\begin{equation*}
\begin{split}
    \ddt{t} {}^R\!\!\Area_t(M)
    &=
    \FP_{z=0} \ddt{t} \int_M x^z \dvol_t
    = \FP_{z=0} \int_M x^z  \omega'(t) \dvol_t \\
    &= \FP_{z=0} \int_M x^z (\sC - R_{g(t)}) \dvol_t \\
    &= \sC \lrpar{ {}^R\!\!\Area_t(M) } - \overset{R\;\;\;}{\int_{M}} R_{g(t)} \dvol_t \\
    &= \sC \lrpar{ {}^R\!\!\Area_t(M) } - 4\pi \chi(M)
\end{split}
\end{equation*}
from which the result follows when
${}^R\!\!\Area_0(M)= 0.$  If $\sC \neq 0,$ then this implies that for some constant $A,$
\begin{equation*}
    {}^R\!\!\Area_t(M) = A e^{\sC t} + \frac{4\pi \chi(M)}{\sC}.
\end{equation*}
We can find $A$ by setting $t=0,$
\begin{equation}\label{VolEq}
    {}^R\!\!\Area_t(M) = \lrspar{ {}^R\!\!\Area_0(M) - \frac{4\pi \chi(M)}{\sC} } e^{\sC t} + \frac{4 \pi \chi(M)}{\sC}
\end{equation}
and the result follows.
\end{proof}

One other natural choice is to take $\sC=r_{i}$ so that the
asymptotic behavior of the curvature at $Y_{i}$ is preserved along
the flow, see Corollary \ref{se.16} below.  This choice is
particularly useful to study the behavior of the metric and the
curvature at infinity.

\subsection{Asymptotic behavior of a solution at infinity} $ $\newline
Given a metric $g_{0}$ satisfying \eqref{ModelFun} and \eqref{ModelCusp}
and $g(t)$ a solution
to \eqref{RicciFlowEq}, we would like to know when these asymptotic behaviors
will be preserved along the normalized Ricci flow.  It turns out to be
convenient to study the asymptotic behavior of the metric $g(t)$ in terms
of the solution $\omega(t)$ of  \eqref{se.1}, since in this setting one
can easily invoke the maximum principle.
The following elementary lemma and proposition are essentially taken from
\cite[\S3.4]{TaoBlog}. We include their proofs for completeness.

\begin{lemma}
Let $t\mapsto h(t)$ be a smooth family of complete Riemannian metrics on
$M$ for $t\in [0,T]$ with curvature uniformly bounded.  Let
$u,v\in \mathcal{C}^{2}([0,T]\times M)\cap \mathcal{C}^{1}([0,T]\times \bar M)$
be two functions such that $u\ge v$ on $[0,T]\times \pa \bar M$
and $u(0,m)\ge v(0,m)$ for all $m\in M.$  Given $A\in \bbR,$ exactly one
of the following two possibilities happens:
\begin{itemize}
\item[(i)] $u(t,m) \ge v(t,m)$ for all $(t,m)\in [0,T]\times \bar M$ or
else
\item[(ii)] there exists $(t,m)\in (0,T]\times M$ such that
\begin{equation*}
\begin{array}{ll}
  u(t,m)< v(t,m), &  (\Delta_{h(t)}u)(t,m)\le (\Delta_{h(t)}v)(t,m), \\
  \nabla u(t,m)= \nabla v(t,m), &
  \frac{du}{dt}(t,m)\le \frac{dv}{dt}(t,m)- A(v(t,m)-u(t,m)).
\end{array}
\end{equation*}
\end{itemize}
\label{se.4}\end{lemma}
\begin{proof}
Replacing $u,v$ with $u-v,0,$ we can assume $v=0.$  Replacing $u$ by
$e^{At}u,$ we can also assume that $A=0.$

In that case, if (i) holds, then clearly (ii) cannot hold.
Conversely,
if $(i)$ does not hold, then the minimum of $u$ is negative and there
exists $(t,m)\in (0,T]\times \bar M$ where this minimum is achieved.  Since we assume that $u\ge v$ on $[0,T]\times \pa \bar M,$ this point is in $(0,T] \times M$ and,
at this point,
we have
\begin{equation*}
\begin{array}{ll}
  u(t,m)< 0, &   \nabla u(t,m)= 0, \\
 (\Delta_{h(t)}u)(t,m)\le 0,&
  \frac{du}{dt}(t,m)\le 0 ,
\end{array}
\end{equation*}
so that $(ii)$ holds.

\end{proof}

\begin{proposition}
With the same assumptions as in Lemma~\ref{se.4}, suppose that $u$ and $v$
also satisfy
\begin{equation*}
\begin{gathered}
\frac{du}{dt}\ge -\Delta_{h(t)}u + \nabla_{X(t)}u + F(t,m,u), \\
\frac{dv}{dt}\le -\Delta_{h(t)}v + \nabla_{X(t)}v + F(t,m,v), \\
\end{gathered}
\end{equation*}
for all $(t,m)\in [0,T]\times M$ where $t\mapsto X(t)$ is a smooth family
of smooth vector fields and $F$ is a function which is uniformly Lipschitz in
the last variable.  Then $u(t,m)\ge v(t,m)$ for all $(t,m)$ in
$[0,T]\times M.$
\label{se.5}\end{proposition}
\begin{proof}
Subtracting the second equation from the first equation and using the
Lipschitz property of $F,$ we get
\begin{equation}
 \frac{d}{dt}(u-v)\ge - \Delta_{h(t)}(u-v) + \nabla_{X(t)}(u-v)
  -C|u-v|
\label{se.6}\end{equation}
where $C$ is the Lipschitz constant of $F.$  Choosing $A>C$ in
Lemma~\ref{se.4}, we see that (ii) cannot occur and hence (i) holds.
\end{proof}

We will first put this into use to study the asymptotic behavior of the
scalar curvature along the flow.

\begin{proposition}
Fix $i\in\{1,\ldots, n_{\hc}+n_{\fun}\}$ and let $\omega$ be a smooth solution to \eqref{se.1} with $\sC=r_{i}= -2e^{-\phi_{i}}$, where we suppose that the
function $\varphi$ of Definition~\ref{se.3a} satisfies \eqref{se.3b}
for some $\delta>0$ and
$k\ge 2.$  If $\omega(t)$ is in $\mathcal{C}^{k+2}_{\Fhc}(M)$
for $t\in [0,T]$, then for all $t\in [0,T]$,
\[
    R_{g(t)}-r_{i}(0) \in x_{i}^{\delta} \mathcal{C}^{k-2}_{\Fhc}(M).
\]
\label{asc.1}\end{proposition}
\begin{proof}
The evolution equation of $R_{g(t)}-r_{i}$ is given by
\[
     \frac{\pa}{\pa t} (R_{g(t)}-r_{i})= -\Delta_{g(t)}(R_{g(t)}-r_{i})
+ R_{g(t)}(R_{g(t)}-r_{i}).
\]

For $\nu>0,$ consider $\psi= x_{i}^{\nu}(R_{g(t)}-r_{i}).$
Then its evolution equation is given by
\[
\frac{\pa \psi}{\pa t} = -\Delta_{g(t)}\psi
+ \nabla_{X(t)}\psi + f\psi, \quad  \left. \psi\right|_{Y_{i}}=0,
\]
where $X^{p}(t)= 2x_{i}^{\nu}g^{pq}\nabla_{q}x_{i}^{-\nu}$ is a
family of vector fields in $\mathcal{C}^{k+2}_{\Fhc}(M;TM)$ and
$f= -x_{i}^{\nu}\Delta_{g(t)}x_{i}^{-\nu} +R_{g(t)}$ is in
$\mathcal{C}^{k}_{\Fhc}(M).$

Since $R_{g(t)}$ is uniformly bounded, we can choose positive constants $C$ and
$C_{1}$ such that $v=Ce^{C_{1}t}x^{\delta+\nu}_i$ satisfies
\begin{equation*}
\begin{gathered}
  \frac{\pa v}{\pa t}= C_{1}v \ge -\Delta_{g(t)}v +\nabla_{X(t)}v+
fv, \\
\frac{\pa (-v)}{\pa t}= -C_{1}v \le
-\Delta_{g(t)}(-v)+\nabla_{X(t)}(-v) +f(-v),
\end{gathered}
\end{equation*}
and $v(t,m)\ge |\psi(t,m)|$ for all $t\in[0,T]$ and $m$ outside a
collar neighborhood of $Y_{i}.$  This last property is to
insure we control what happens at the other ends of the surface.
  Choosing $C>0$ big enough,
we can also assume that $v(0,m)\ge |\psi(0,m)|$ for all
$m\in M.$  We can then apply Proposition~\ref{se.5}
to conclude that
\[
     -v(t,m)\le \psi(t,m) \le v(t,m)
\]
for all $(t,m)\in [0,T]\times M.$  Thus, this gives that
\[
    -C e^{C_{1}t}x_{i}^{\delta} \le R_{g(t)}-r_{i}\le
Ce^{C_{1}t}x_{i}^{\delta}.
\]
We can derive similar estimates for the derivatives of $R_{g(t)}$ (up to
order $k-2$) by looking at their evolution equations, from which the result follows.
\end{proof}

More generally, unless we have that $r_{i}=\sC$ to start with,
the asymptotic value $r_{i}$ of the scalar curvature in the end $Y_{i}$ will vary with $t$.  To emphasize this potential dependence on $t$, we will use the notation $r_{i}(t).$
It is easy to see that for $\sC<0$
\begin{equation}
    r_i'(t)  = - r_i(t) ( \sC - r_i(t))
    \implies
    \frac{ \abs{ \sC - r_i(t) } } { \abs{r_i(t)} }
    = \frac{ \abs{ \sC - r_i(0) } } { \abs{r_i(0)} } e^{\sC t} \notag
\label{FirstLine} \end{equation}
which
implies that
\begin{equation}\label{FormulaAsymCurv}
\begin{gathered}
    \frac{ \sC - r_i(t)  } { r_i(t) }
    = \frac{ \sC - r_i(0) } { r_i(0) } e^{\sC t}
    \implies
    r_i(t) = \frac{ r_i(0) \sC } { r_i(0) + (\sC - r_i(0)) e^{\sC t} }
\end{gathered}
\end{equation}

\begin{proposition}\label{PreserveDecay}
  Let $\omega$ be a smooth solution to \eqref{se.1} with $\varphi\in
\mathcal{C}^{k+2}_{\Fhc}(M)$ for some $k\ge 2$ satisfying \eqref{se.3b}
 for some $\delta>0.$  Suppose
$\omega(t)$ is in $\mathcal{C}^{k+2}_{\Fhc}(M)$ for $t\in
[0,T].$ Then for each $i,$ there exists a smooth function
$c_{i}:[0,T]\to \bbR$ such that
\begin{equation}
   \omega(t)-c_{i}(t)\in x^{\delta}_{i}\mathcal{C}^{k-2}_{\Fhc}(M)\; \forall \
t \in [0,T].
\label{se.11}\end{equation}
\end{proposition}

\begin{proof}
Fix $i\in\{1,\ldots,n_{\hc}+n_{\fun}\}.$  Recall that under the rescaling in time and space given by
\[
   \tau= \frac{e^{-\sC t}-1}{-\sC}, \quad h= e^{-\sC t} g,
\]
the normalized Ricci flow equation $\frac{\pa g}{\pa t}= (\sC-R_{g})g$ becomes the usual Ricci flow,
\begin{equation}
     \frac{\pa h}{\pa \tau} = -R_{h} h.
\label{usualRicci}\end{equation}
Similarly, under the change of variable $\tau= \frac{e^{-\sC' t'}-1}{-\sC'}$, $h= e^{-\sC't'}g'$, where $\sC'$ is a non-zero constant,   equation~\eqref{usualRicci} becomes
\begin{equation*}
       \frac{\pa g'}{\pa t'} = (\sC'-R_{g'})g'(t').
\end{equation*}
This means that at the cost of  rescaling in time and space if necessary,
we are free to choose the constant $\sC$ as we want in \eqref{omegaEq}.  To study the behavior of $\omega$ near $Y_{i}$, the best choice is to take $\sC=r_{i}(0)$.     In that case, we need to show that
$\omega(t)\in x_{i}^{\delta}\mathcal{C}^{k-2}_{\Fhc}(M)$ for all
$t\in [0,T].$
Let $\nu\in (0,\delta)$ be given.  Then the function $\psi_{i}= x^{\nu}_{i}\omega$
satisfies the equation
\begin{equation}
\begin{gathered}
  \frac{\pa \psi_{i}}{\pa t} = -e^{-\omega}(\Delta_{g_{0}}\psi_{i} +
x^{\nu}_{i}R_{g_{0}}) + x^{\nu}_{i}\sC + e^{-\omega}f \psi_{i} +
\nabla_{X(t)}\psi_{i},  \\
\left. \psi\right|_{\pa \bar M}=0,
\end{gathered}
\label{se.12}\end{equation}
where $X(t)\in \mathcal{C}^{k+2}_{\Fhc}(M;TM)$ is a family of vector fields on $M$ and $f\in \mathcal{C}^{k}_{\Fhc}(M).$
Notice in particular that
$\sup_{M}|X(t)|_{g(0)}<\infty$ for all $t\in [0,T].$
  We can rewrite this equation
as
\begin{equation}
  \frac{\pa \psi_{i}}{\pa t}= -\Delta_{g(t)}\psi_{i} + \nabla_{X(t)}\psi_{i}
+ f_{1}(t,m)\psi_{i} + f_{0}(t,m)
\label{se.13}\end{equation}
with $f_{0}(t,m)= -e^{-\omega(t,m)}x^{\nu}_{i}R_{g_{0}}+ x^{\nu}_{i}\sC
$ and $f_{1}(t,m)=e^{-\omega(t,m)}f.$
We know by Proposition~\ref{asc.1} and \eqref{se.7} that
\begin{equation}
  f_{0}(t,\cdot)\in x^{\nu+\delta}\mathcal{C}^{\infty}_{\Fhc}(M)
  \quad \forall \; t\in [0,T].
\label{se.14}\end{equation}
Thus, since $\omega$ is uniformly bounded, we can choose positive constants
$C,C_{1}$  sufficiently
large so that $v=Ce^{C_{1}t}x^{\delta+\nu}$ satisfies
\begin{equation}
\begin{gathered}
  \frac{\pa v}{\pa t}= C_{1}v \ge -\Delta_{g(t)}v +\nabla_{X(t)}v+
f_{1}v +f_{0}, \\
\frac{\pa (-v)}{\pa t}= -C_{1}v \le -\Delta_{g(t)}(-v)+\nabla_{X(t)}(-v) +f_{1}(-v)+ f_{0},
\end{gathered}
\label{se.15}\end{equation}
and $v(t,m)\ge |\psi_{i}(t,m)|$ for $t\in [0,T]$ and $m$ outside a collar
neighborhood of $Y_{i}$ in $\bar M.$  We can also choose
$C>0$ so that $v(0,m)\ge |\psi_{i}(0,m)|$ for all $m\in M.$
We can then apply
Proposition~\ref{se.5} to conclude that
\[
             -v(t,m)\le \psi_{i}(t,m) \le v(t,m), \quad \forall (t,m)\in [0,T]\times M.
\]
Thus, this gives that
\[
             -Ce^{C_{1}t}x^{\delta}\le \omega \le Ce^{C_{1}t}x^{\delta}.
\]
We can derive similar estimates for the derivatives of $\omega$ (up to order
k-2) by looking
at their evolution equations (which can be derived by using the
identity $\nabla \Delta= \Delta \nabla + \frac{1}{2}R\nabla$), from which the result follows.
\end{proof}

\begin{corollary}
Let $\omega(t)$ be as in Proposition~\ref{PreserveDecay}. Then the scalar
curvature $R_{g(t)}$ of the metric $g(t)=e^{\omega(t)}g_{0}$ is such that
\[
   R_{g(t)}- r_{i}(t)\in x^{\delta}_{i}\mathcal{C}^{k}_{\Fhc}(M) \quad
   \forall \; t\in[0,T],
\]
where $r_{i}:[0,T]\to \bbR$ is given by \eqref{FormulaAsymCurv}.
When $\sC=r_{i}(0),$ then  $r_{i}$ is constant along the flow.
\label{se.16}\end{corollary}
\begin{proof}
This is a direct consequence of Proposition~\ref{asc.1}
and Proposition~\ref{PreserveDecay}.
\end{proof}

When we will apply Polyakov's formula to a family of metrics $g(t)=e^{\omega(t)}g_{0}$ evolving according to the Ricci flow, it will be convenient to know
that the conformal factor $\omega(t)$ remains smooth up to the boundary
along the flow.  This is the content of the next proposition.

\begin{proposition}\label{PreserveSmoothness}
Let $\omega$ be a smooth solution to
\eqref{se.1} with initial metric $g_{0}$ satisfying
\eqref{ModelFun} and \eqref{ModelCusp} with $\varphi\in \CI_{\Fhc}(M)\cap \CI(\bar M)$ satisfying \eqref{se.3b} for some $\delta>0.$
Suppose that $\omega(t)$ is uniformly in
$\CI_{\Fhc}(M)$ for $t\in [0,T].$  Then
\[
            \omega(t)\in \CI(\bar M)\cap \CI_{\Fhc}(M)\;\forall t\in [0,T].
\]
\label{bph.1}\end{proposition}

\begin{proof}
Notice that if there are no cusp ends, then
$\CI(\bar M)\subset \CI_{\Fhc}(M).$  However, if
$\pa_{\hc}\bar M\ne \emptyset,$ then this inclusion is false.  In fact,
near a cusp end, we have that
\begin{equation}
  f\in \CI_{\Fhc}(M)\cap \CI(\bar M) \quad \Longrightarrow
\left. f\right|_{\pa_{\hc}\bar M} \; \mbox{is locally constant.}
\label{bph.2}\end{equation}
To show that $\omega(t)\in \CI(\bar M),$ we will inductively construct the
Taylor series of $\omega$ at the boundary.  Thus, we need to show that
there exists $\omega_{k}\in \CI([0,T]\times \pa \bar M)$ for $k\in \bbN\cup\{0\}$
such that
\begin{equation}
  \left(\omega(t)- \sum_{k=0}^{N} \chi(x)\omega_{k}(t)x^{k}\right)  \in x^{N+1}\CI_{\Fhc}(M)
\quad \forall \; N\in \bbN\cup\{0\}
\label{bph.3}\end{equation}
in a collar neighborhood of $\pa \bar M,$ where $\chi:[0,+\infty) \to
[0,+\infty)$ is a smooth cut-off function with
$\chi(x)=1$ for $x<\frac{\epsilon}{2}$ and $\chi(x)=0$ for
$x>\frac{3\epsilon}{4}.$

If $\omega_{k}$ exists, notice by \eqref{bph.2} that it
is locally constant on $\pa_{\hc} \bar M.$
Without loss of generality, we can assume that $\pa \bar M$ has only one boundary component which
is associated either to a cusp or a funnel.
We first need to define $\omega_{0}.$  If $\omega(t)$ were in $\CI(\bar M)$
as claimed, then the evolution equation for $\omega_{0}$ would be
\begin{equation}
\frac{\pa \omega_{0}}{\pa t}= -e^{-\omega_{0}}r_{i}(t)+\sC, \quad
\omega_{0}(0)\equiv 0.
\label{bph.4}\end{equation}
Thus, we can define $\omega_{0}$ to be the unique solution to
\eqref{bph.4}.  Rescaling the flow and the metric if necessary, we can assume that $\sC=r_{i}=-2$
so that $R_{g(t)}+2 \in x^{\delta} \CI_{\Fhc}(M).$  Then $\tomega_{1}(t)= \omega(t)-\omega_{0}(t)\chi(x)$  satisfies the
evolution equation
\begin{equation}
\begin{aligned}
\frac{\pa \tomega_{1}}{\pa t}&= -e^{-\omega}(\Delta_{g_{0}}\omega +R_{g_{0}})
  + \sC + \chi(x)\left(e^{-\omega_{0}}(-2) -\sC\right) \\
&= -e^{-\omega}\Delta_{g_{0}}\omega + R_{g_{0}}(-e^{-\omega}+
  e^{-\chi(x)\omega_{0}})   \\
& + \left(-e^{-\chi(x)\omega_{0}}R_{g_{0}}-2
  \chi(x)e^{-\omega_{0}}\right)
+\sC(1-\chi(x)) \\
&= -\Delta_{g(t)}\tomega_{1} + \tilde{f}_{1}\tomega_{1}+
\tilde{h}_{1}
\end{aligned}
\label{bph.5}\end{equation}
where $\tilde{f}_{1}=R_{g_{0}}e^{-\chi(x)\omega_{0}}\left(
\frac{1-e^{-\tomega_{1}}}{\tomega_{1}} \right)$ is in $\CI_{\Fhc}(M)$ for all $t\in [0,T]$ and has
the same regularity as $\omega$ at the boundary, while
$\tilde{h}_{1}\in x(\CI_{\Fhc}(M)\cap \CI(\bar M)).$
Thus using a barrier function of the form
$v= Ce^{C_{1}t}x$ for some $C,C_{1}>0$ large enough, we can proceed
as in the proof of Proposition~\ref{PreserveDecay}, to show that
$\bomega_{1}= \frac{\tomega_{1}}{x} \in \CI_{\Fhc}(M)$ for all
$t\in[0,T].$  From \eqref{bph.5}, its evolution equation is of the form
\begin{equation}
  \frac{\pa \bomega_{1}}{\pa t}= -\Delta_{g(t)}\bomega_{1}+
\nabla_{X_{1}(t)}\bomega_{1} + \bar{f}_{1}\bomega_{1} +\bar{h}_{1}(t,m),
\quad \bomega_{1}(0)\equiv 0,
\label{bph.6}\end{equation}
where $X_{1}(t)\in \CI_{\Fhc}(M;TM)$ and
$\bar{f}_{1}(t,\cdot)$ is in $\CI_{\Fhc}(M)$ for
all $t\in [0,T]$ and has the same regularity as $\omega$ at the boundary,
while
$\bar{h}_{1}\in \CI_{\Fhc}(M)\cap\CI(\bar M)$ for all $t\in [0,T].$

Suppose now for a proof by induction that $\omega_{0},\ldots, \omega_{N-1}$
have been defined to satisfy \eqref{bph.3} and that
\begin{equation}
   \bomega_{N}= \frac{\omega -\sum_{k=1}^{N-1}\chi(x)\omega_{k}x^{k}}
{x^{N}}\in \CI_{\Fhc}(M)
\label{bph.7}\end{equation}
satisfies the evolution equation
\begin{equation}
 \frac{\pa \bomega_{N}}{\pa t}= -\Delta_{g(t)} \bomega_{N}
+ \nabla_{X_{N}(t)}\bomega_{N}+
\bar{f}_{N}\bar \omega_{N}+ \bar{h}_{N}
\label{bph.8}\end{equation}
where $X_{N}(t)\in \CI_{\Fhc}(M;TM)$ and with $\bar{f}_{N}(t,\cdot),$
$\bar{h}_{N}(t,\cdot)$ in $\CI_{\Fhc}(M)$ having the
same regularity as $\bomega_{N-1}$ at the boundary (with the convention
that $\bomega_{0}=\omega$).  We then define
$\omega_{N}\in\CI([0,T]\times \pa\bar M)$ to be the unique solution
of the evolution equation
\begin{equation}
  \frac{\pa \omega_{N}}{\pa t} = \left(\left.
\bar{f}_{N}\right|_{\pa \bar M}\right) \omega_{N} +
\left. \bar{h}_{N}\right|_{\pa \bar M}, \quad \omega_{N}(0)\equiv 0.
\label{bph.9}\end{equation}
Then $\tomega_{N+1}= \bomega_{N}-\chi(x)\omega_{N}$ satisfies the
evolution equation
\begin{equation}
\begin{aligned}
 \frac{\pa \tomega_{N+1}}{\pa t}&= -\Delta_{g(t)}\bomega_{N}+ \nabla_{X_{N}(t)}\bomega_{N}+\bar{f}_{N}\bomega_{N}+ \bar{h}_{N} \\
& -\chi(x)\left(\left.\bar{f}_{N}\right|_{\pa\bar M}\right)\omega_{N}
-\chi(x)\left( \left. \bar{h}_{N}\right|_{\pa \bar M}\right) \\
&= -\Delta_{g(t)}\bomega_{N} +
\nabla_{X_{N}(t)}\bomega_{N} +
\bar{f}_{N}(\bomega_{N}-\chi(x)\omega_{N}) \\
&+
(\bar{f}_{N}-\left.\bar{f}_{N}\right|_{\pa \bar M})\chi(x)\omega_{N}
+ \left(\bar{h}_{N} - \chi(x)\left.\bar{h}_{N}\right|_{\pa \bar M}\right)
\\
&= -\Delta_{g(t)}\tomega_{N+1}+ \nabla_{X_{N}(t)}\tomega_{N+1}
+ \bar{f}_{N}\tomega_{N+1}+ \tilde{h}_{N+1}
\end{aligned}
\label{bph.10}\end{equation}
where
$\tilde{h}_{N+1}(t,\cdot)\in x\CI_{\Fhc}(M)$ has
the same regularity as $\bomega_{N-1}$ at the boundary.
Thus, using a barrier function of the form $v= Ce^{C_{1} t}x$ for
some $C,C_{1}>0$ large enough, we can again proceed as in the proof of
Proposition~\ref{PreserveDecay} to show that
$\bomega_{N+1}(t):= \frac{\tomega_{N+1}}{x}$ is in
$\CI_{\Fhc}(M).$  Moreover, it satisfies
the evolution equation
\begin{equation}
\frac{\pa \bomega_{N+1}}{\pa t}= \Delta_{g(t)}\bomega_{N+1} +
\nabla_{X_{N+1}(t)}\bomega_{N+1}+ \bar{f}_{N+1}\bomega_{N+1}
+\bar{h}_{N+1}
\label{bph.11}\end{equation}
with $\bar{f}_{N+1}(t,\cdot),\bar{h}_{N+1}(t,\cdot) \in \CI_{\Fhc}(M)$
having the same
regularity up to the boundary as $\bomega_{N}$ and
with $X_{N+1}(t)\in \CI_{\Fhc}(M;TM).$

In this way, we define inductively $\omega_{k}\in
\CI([0,T]\times\pa \bar M)$ such
that \eqref{bph.3} is satisfied for all $N\in \bbN.$

\end{proof}

\subsection{Short-time existence and uniqueness}

The short-time existence of a solution of \eqref{se.1} follows from the more
general result of Shi \cite{Shi1989} who established the short-time existence
of the Ricci flow on a complete Riemannian manifold with bounded curvature.
In our context, since the equation can be written in the scalar form
\eqref{se.1}, it is possible to prove uniqueness relatively easily using
the maximum principle.

\begin{proposition}[Short-time existence and uniqueness]
Consider the equation \eqref{se.1} with initial metric $g_{0}$ satisfying
\eqref{ModelFun} and \eqref{ModelCusp} with
$\varphi\in \mathcal{C}^{k+2}_{\Fhc}(M)$
satisfying \eqref{se.3b} for some $\delta>0$ and $k\ge 2.$  In particular, there exists $K>0$ such that $|R_{g_{0}}|<K$
everywhere on $M.$  Then there exists $T=T(K,\sC)>0$ depending
on $K$ and $\sC$ such that \eqref{se.1} has a unique smooth
solution $\omega(t)$ on $[0,T]\times M$ satisfying
\[
    \omega(t)-c_{i}(t)\in x^{\delta}_{i}\mathcal{C}^{k-2}_{\Fhc}(M)
\]
for some smooth functions $c_{i}:[0,T]\to \bbR.$
\label{se.17}\end{proposition}
\begin{proof}
Short-time existence follows from the short-time existence result of
Shi \cite{Shi1989} for the Ricci flow on complete manifolds with bounded
curvature.  The decay property of the solution follows from
Proposition~\ref{PreserveDecay}.

To prove uniqueness, suppose $\omega_{1}$ and $\omega_{2}$ are two solutions.
By the proof of Proposition~\ref{PreserveDecay} (see also \eqref{bph.4}), we
know that $\left.\omega_{1}\right|_{\pa \bar M}= \left. \omega_{2}\right|_{\pa \bar M}.$
In fact, by Proposition~\ref{PreserveDecay}, we have that $v=\omega_{1}-\omega_{2}$
is in $x^{\delta}\mathcal{C}^{k-2}_{\Fhc}(M)$ and satisfies
\begin{equation}
\left\{ \begin{array}{l}
\frac{\pa v}{\pa t} = -\Delta_{g_{1}(t)}v
 +F(t,m,\omega_{1})-F(t,m,\omega_{2}),  \quad g_{1}(t)= e^{\omega_{1}(t)}g_{0}, \\
  v(0)\equiv 0,
\end{array}     \right.
\label{se.18}\end{equation}
where
\begin{equation}
  F(t,m,\psi)= (-\Delta_{g_{0}}\omega_{2}(m,t)-R_{g_{0}}(m))e^{-\psi},
\quad t\in [0,T], \; m\in M, \; \psi\in \bbR.
\label{se.19}\end{equation}
Since $\omega_{1}$ and $\omega_{2}$ together with their  derivatives are
uniformly bounded on $[0,T],$ there exists a constant $A>0$
 such that
\begin{equation}
  |F(t,m,\omega_{1})- F(t,m,\omega_{2})| \le A |v(t,m)| \quad \forall \;
          (t,m)\in [0,T]\times M.
\label{se.20}\end{equation}
Consequently,
\begin{equation}
  \frac{\pa v}{\pa t} \ge -\Delta_{g_{1}(t)}v - A|v|, \quad v(0)\equiv 0,
\quad \left. v\right|_{\pa M}=0.
\label{se.21}\end{equation}
By Proposition~\ref{se.5} applied to $v$ and $0,$
this means  that $v(t,m)\ge 0$ for all
$(t,m)\in [0,T]\times M.$  Interchanging the r\^oles of $\omega_{1}$ and
$\omega_{2},$ we can also show that $v\le 0,$ which implies that
$v\equiv 0,$ that is, $\omega_{1}\equiv \omega_{2}$ on
$[0,T]\times M,$ establishing uniqueness.
\end{proof}

\subsection{Long-time existence}

To prove long-time existence, it suffices to get an a priori bound on the
scalar curvature, for then we can apply Proposition~\ref{se.17} recursively
to get long-time existence.  As noticed by Hamilton \cite{Hamilton},
on a compact surface, such an a priori estimate follows from the existence
of a potential function.  This approach was recently generalized
to surfaces with cusps ends by Ji, Mazzeo and
\v{S}e\v{s}um \cite{Ji-Mazzeo-Sesum}.

To get long-time existence in our case, we need to generalize the
construction of the potential function given in \cite{Ji-Mazzeo-Sesum} to
also include funnel ends.

\begin{proposition}[Potential function]
Let $g$ be a Riemannian metric on $M$ satisfying \eqref{ModelFun} and
\eqref{ModelCusp}
with $\varphi\in \mathcal{C}^{3}_{\Fhc}(M)$ satisfying
\eqref{se.3b} for
some $\delta_{\fun}>2$ and $\delta_{\hc}>0.$
  Let $\sC<0$ be arbitrary if
$\pa_{\fun}M\ne \emptyset$ and $\sC= \bar R$ if
$\pa_{\fun}\bar M =\emptyset.$
Then there exists a unique function
$f$ and constants $c_{i}$ such that
\begin{equation}
  -\Delta_{g}f= R -\sC, \quad (f- \sum_{i=1}^{n_{\fun}+n_{\hc}}c_{i}\log x_{i})\in
C^{2,\alpha}_{\Fhc}(M),
  \quad \sup_{M} |\nabla f|_{g}<\infty,
\label{le.1a}\end{equation}
with $\int_{M}fdg=0$  if $\pa_{\fun}\bar M=\emptyset$ and
\[
\left.(f- c_{1}\log x_{1})\right|_{Y_{1}}=0,
\quad \left. \frac{\pa}{\pa x_{j}}(f- \sum_{i=2}^{n_{\fun}}c_{i}\log x_{i})\right|_{Y_{j}}=0, \quad j\in\{2,\ldots,n_{\fun}\}
\]
otherwise.

\label{le.1}\end{proposition}

\begin{proof}
If $M$ is compact, this is very easy to prove.  If $(M,g)$ has only cusp
ends, this was proved in \cite{Ji-Mazzeo-Sesum}, in fact with weaker decay
assumptions on $\varphi.$  To prove the proposition when there are funnel ends,
the idea is to reduce to the case where there are only cusp ends (or when
$M$ is compact) via a doubling construction.

Thus let us start with a $\Fhc$-metric $g$ for which $\pa_{\fun}\bar M\ne \emptyset.$
\begin{figure}[h]
\centerline{\includegraphics[height=5cm]{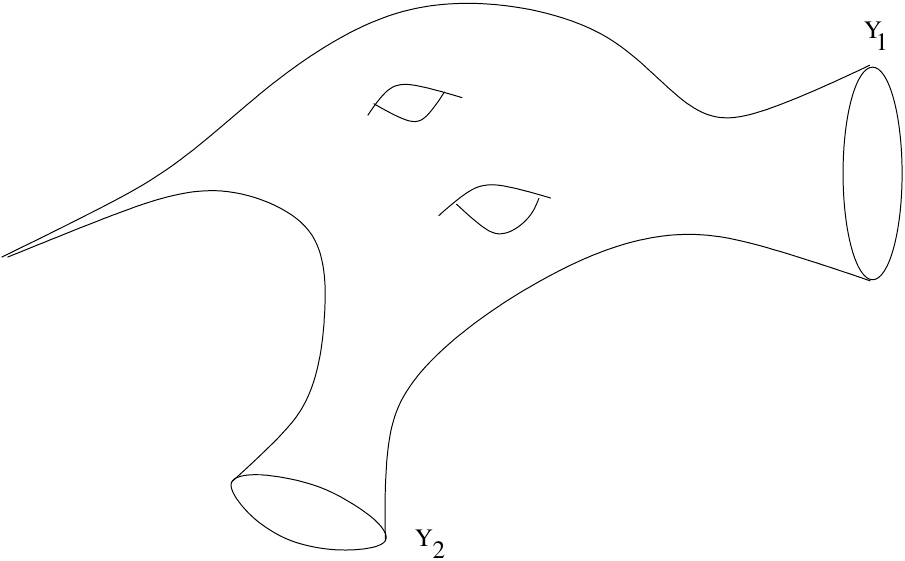}}
\caption{$(M,g)$} \label{fig1}
\end{figure}

By \eqref{se.7}, we know that for $i\in\{1,\ldots, n_{\fun}\},$
\begin{equation}
    R_{g}-r_{i} \in x_{i}^{\delta_{\fun}}C^{1}_{\Fhc}(M).
\label{pot.1}\end{equation}
By assumption, there also exists a constant $\phi_{i}$ such that
\begin{equation}
    (\varphi-\phi_{i})\in x_{i}^{\delta_{\fun}}C^{3}_{\Fhc}(M).
\label{pot.2}\end{equation}
Then the function
\begin{equation}
\psi_{i}=e^{\phi_{i}}(r_{i}-\sC)\log x_{i}
\label{pot.2b}\end{equation}
is such that
\begin{equation}
(-\Delta \psi_{i}-(R-\sC))\in x_{i}^{\delta_{\fun}}C^{1}_{\Fhc}(M).
\label{pot.3}\end{equation}
Moreover, since we chose $x_{i}$ to be equal to one outside a collar neighborhood of $Y_{i}$, we see that $\psi_{i}$ has its support contained in a collar neighborhood of $Y_{i}$.

Thus, if we set $\tilde{f}= f- \sum_{i=1}^{n_{\fun}}\psi_{i},$ we
can rewrite $-\Delta f=R-\sC$ as
\begin{equation}
      -\Delta \tilde{f}= h, \quad \mbox{with} \; h=
(R-\sC + \sum_{i=1}^{n_{\fun}}\Delta \psi_{i}) \in
\mathcal{C}^{1}_{\Fhc}(M).
\label{pot.4}\end{equation}

In a collar neighborhood of $\pa_{\fun}\bar M,$ equation \eqref{pot.4}
 takes the form
\begin{equation}
  -e^{-\varphi} x^{2}_{\fun}\Delta_{g_{E}}\tilde{f}= h
\label{le.2}\end{equation}
where $\Delta_{g_{E}}$ is the Laplacian associated to the incomplete
cylindrical metric $g_{E}= dx^{2}+ \pi_{F}^{*}h_{F}$ where $h_{F}$ is
a metric on $\pa_{\fun}\bar M$ and $\pi_{F}$ is the projection from the
collar neighborhood of $\pa_{\fun}\bar M$ onto $\pa_{\fun}\bar M$.

Thus, with respect to the
metric
\begin{equation}
   \tilde{g}= e^{-\varphi\chi}x^{2}_{\fun}g,
\label{le.5}\end{equation}
where $\chi\in \CI_{c}(\pa_{F}M\times [0,\epsilon)_{x_{\fun}})$ is a
nonnegative
cut-off function equal to 1 for $x_{\fun}<\frac{\epsilon}{2}$ and
equal to zero for $x_{\fun}> \frac{3\epsilon}{4},$
 equation \eqref{le.2} can be rewritten as
\begin{equation}
  -\Delta_{\tilde{g}}\tilde{f}= \tilde{h}, \quad \tilde{h}=  \frac{h}{e^{-\varphi\chi} x^{2}_{F}}.
\label{le.6}\end{equation}

\begin{figure}[h]
\centerline{\includegraphics[scale=.7]{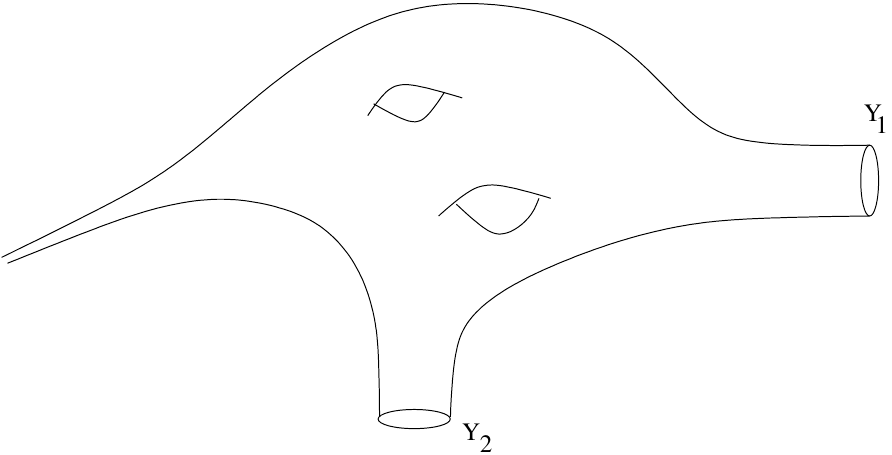}} 
\caption{$(M,\tilde{g})$} \label{fig2}
\end{figure}

The metric $\tilde{g}$ is incomplete and if we glue two copies
of $\bar M$ along $Y_{1}$ to get
\begin{equation}
  \bbM_{1}= \bar M \cup_{Y_{1}} \bar M,
\label{le.8}\end{equation}
then the metrics $\tilde{g}$ on each copy of $\bar M$ glue together to
give a smooth metric $\hat{g}_{1}$ on $\bbM_{1}.$
The Riemannian manifold $(\bbM_{1},\hat{g}_{1})$ has $2(n_{\fun}-1)$
ends where the metric is asymptotic to an incomplete cylinder and
$2n_{\hc}$ ends were it is asymptotic to a cusp.

\begin{figure}[h]
\centerline{\includegraphics[width=13.5cm]{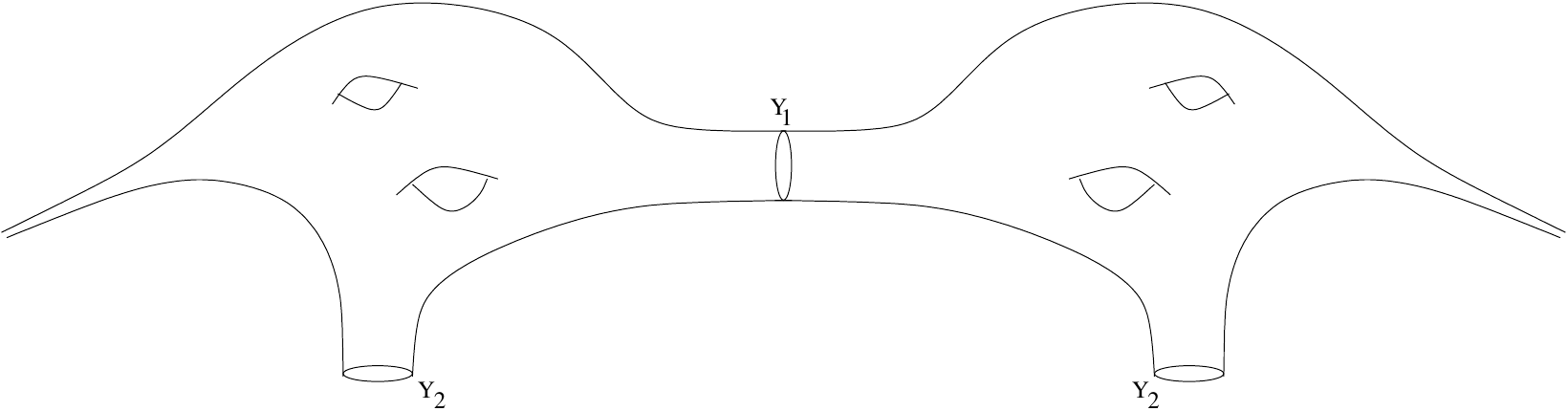}} 
\caption{$(\bbM_{1},\hat{g}_{1})$} \label{fig3}
\end{figure}

Let
\begin{equation}
  \pa_{F}\bbM_{1}\cong \left(\bigcup_{i=2}^{n_{\fun}} Y_i \right)\bigsqcup
\left(\bigcup_{i=2}^{n_{\fun}} Y_i \right)
\label{DN.1}\end{equation}
be the part of the boundary associated to cylindrical ends.  We can consider
the double of $\bbM_{1}$ along $\pa_{\fun}\bbM_{1}$
\begin{equation}
     \bbM_{2}= \bbM_{1}\cup_{\pa_{\fun}\bbM_{1}}\bbM_{1}.
\label{DN.2}\end{equation}
The metric $\hat{g}_{1}$ on each copy of $\bbM_{1}$ glue together to give
a smooth metric $\hat{g}_{2}$ on $\bbM_{2}.$
Clearly, $(\bbM_{2},\hat{g}_{2})$ is
a complete surface with $4 n_{\hc}$ cusp ends (or is a compact surface if
$n_{hc}=0$).

\begin{figure}[h]
\centerline{\includegraphics[width=13.5cm]{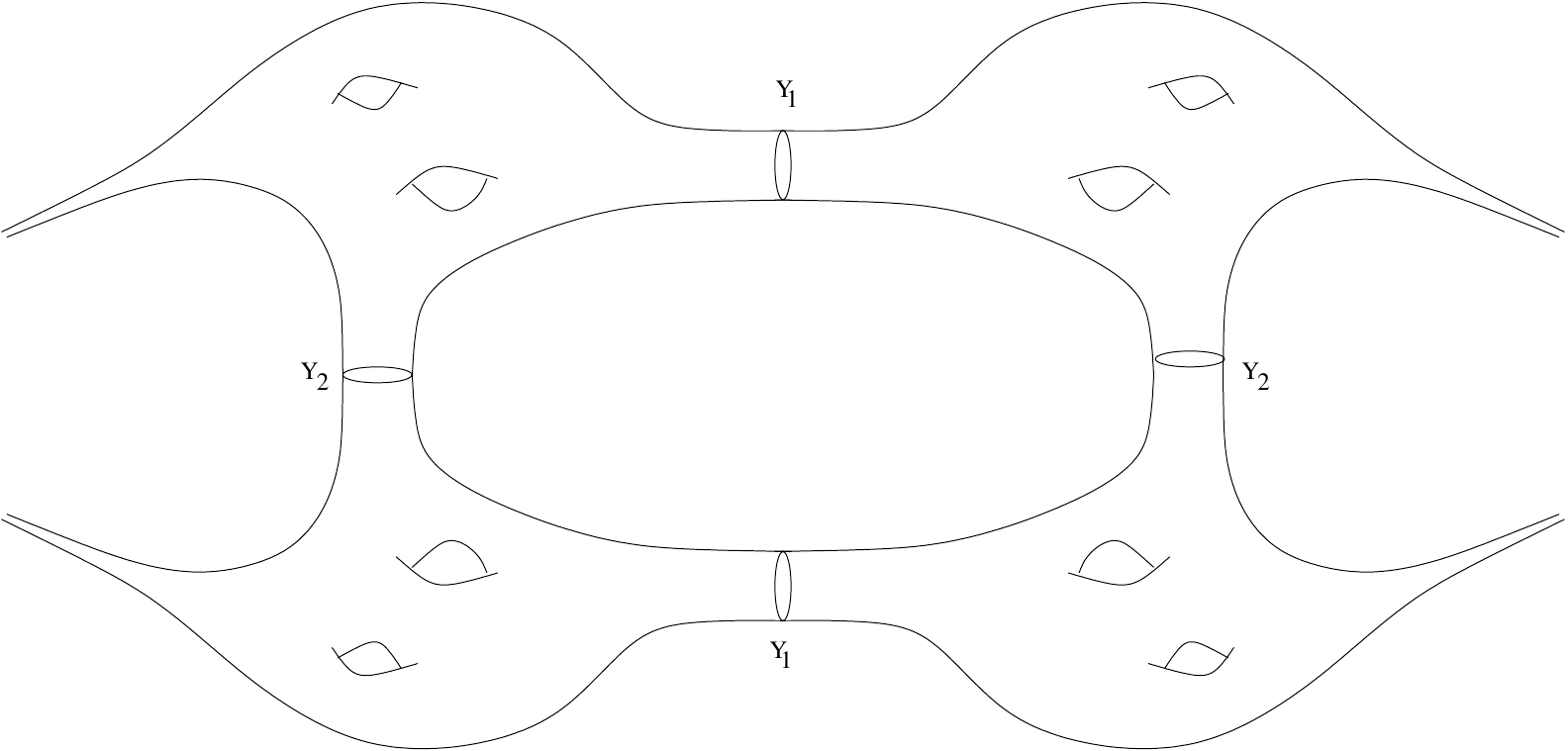}} 
\caption{$(\bbM_{2},\hat{g}_{2})$} \label{fig4}
\end{figure}

Let $\hat{x}\in \CI(\bar\bbM_{2})$ be a boundary defining function
on $\bar\bbM_{2}$ whose restriction to each copy of $\bar M$ in $\bar\bbM_{2}$
is equal to $x_{\hc}.$
Let $\hat{h}_{1}$ be the function on $\bbM_{1}$ whose restriction to
one copy of $M$ in $\bbM_{1}$ is $\tilde{h}$ and whose restriction
to the other copy is $-\tilde{h}.$  Let $\hat{h}_{2}$ be the function whose restriction
to each copy of $\bbM_{1}$ in $\bbM_{2}$ is $\hat{h}_{1}.$   By \eqref{pot.1} and since
$\delta_{\fun}>2,$ we have that $\hat{h}_{2}\in \hat{x}^{\delta_{\hc}}\mathcal{C}^{0,\alpha}_{\hat{g}_{2}}(\bbM_{2}).$ Then on
$\bbM_{2},$ one can consider the equation
\begin{equation}
  -\Delta_{\hat{g}_{2}}\hat{f}= \hat{h}_{2}  \quad \mbox{on}\; \bbM_{2}.
\label{DN.3}\end{equation}
Clearly, by symmetry, $\int_{\bbM_{2}}\hat{h}_{2}d\hat{g}_{2}=0$ so that we
can
apply the result of \cite{Ji-Mazzeo-Sesum} to conclude that there
exists a unique solution $\hat{f}$ with constants $\hat{c}_{i},$
$i\in\{1,\ldots,2n_{\hc}\}$ such that
\[
(\hat{f}- \sum_{i=1}^{4n_{\hc}}\hat{c}_{i}\log{\hat{x}_{i}}) \in
\mathcal{C}^{2,\alpha}_{\hat{g}}(\bbM)
\]
and
\begin{equation}
 -\Delta_{\hat{g}_{2}}\hat{f}= \hat{h}_{2},
\quad \int_{\bbM_{2}}\hat{f}d\hat{g}_{2}=0,
\quad \sup_{\bbM_{2}}|\nabla \hat{f}|_{\hat{g}_{2}}<\infty.
\label{le.10}\end{equation}

Since this
solution is unique, we see by symmetry that the restriction of this
solution to one of the copies of $\bar M$ in $\bbM_{2}$ will solve the
equation $-\Delta_{g}\tilde{f}= R-\sC$ with bounded gradient,
Dirichlet boundary condition on
$Y_{1}$ and Neumann boundary condition on
$\cup_{i=2}^{n_{\fun}}Y_{i}.$
Thus $f= \tilde{f}+ \sum_{i=1}^{n_{\fun}}\psi_{i}$ will be
the desired solution and is clearly unique.

\end{proof}

It is also interesting to consider the following variant of the construction
which only involves Neumann boundary conditions.

\begin{proposition}
Let $g$ be a Riemannian metric satisfying
\eqref{ModelFun} and \eqref{ModelCusp} with
$\varphi\in \CI(\bar M)\cap \cC^3_{\Fhc}(M)$ satisfying
\eqref{se.3b} for $\delta_{\fun}=2$ and some $\delta_{\hc}>0.$  Suppose
that $\pa_{\fun}M\ne 0$ and that ${}^{R}\!\!\Area(g)\ne 0$ so that
$\bar R$ is well-defined.  Then
there exists a unique $f$ and constants $c_{i}$ such that
\begin{equation*}
  -\Delta_{g}f= R -\bar{R}, \quad (f- \sum_{i=1}^{n_{\fun}+n_{\hc}}c_{i}\log x_{i})\in
C^{2,\alpha}_{\Fhc}(M),
  \quad \sup_{M} |\nabla f|_{g}<\infty,
\end{equation*}
with
\[
\left. \frac{\pa}{\pa x_{\fun}}(f- \sum_{i=2}^{n_{\fun}}c_{i}\log x_{i})\right|_{\pa_{\fun}\bar M}=0 \quad \mbox{and}  \quad \int_{M}x_{\fun}^{2}f dg =0.
\]
\label{DN.4}\end{proposition}
\begin{proof}
We proceed as in the proof of Proposition~\ref{le.1} to define the functions
$\psi_{i},$ $\tilde{h}$ and $\tilde{f}$ and the metric $\tilde{g}.$  However,
instead of $\bbM_{1}$ and $\bbM_{2},$ we consider directly the double of
$M$ along $\pa_{\fun}\bar M,$
\begin{equation*}
     \bbM = \bar M \cup_{\pa_{\fun}\bar M}\bar M.
\end{equation*}
The metrics $\tilde{g}$ on each copy of $M$ glue together to give a
metric $\hat{g}$ on $\bbM$ having only cusp ends.
On $\bbM,$ we consider the function $\hat{h}$ whose restriction to each
copy of $\bar M$ in $\bbM$ is $\tilde{h}$ and the equation
\begin{equation}
   -\Delta_{\hat{g}}\hat f= \hat{h}.
\label{DN.6}\end{equation}
Since $\varphi\in \CI(\bar M),$ we only need $\delta_{\fun}=2$ to deduce from
\eqref{pot.1} that $\hat{h}\in \hat{x}^{\delta_{\hc}}\mathcal{C}^{1}_{\hat{g}}(\bbM).$
A quick computation shows that $\overset{R\;\;\;\;}{\int_{M}} \Delta_{g}\psi_{i}dg=0$ so that
\begin{equation}
\int_{\bbM}\hat{h}d\hat{g}= 2\int_{M}\tilde{h}d\tilde{g}=
2\int_{M}h dg= 2 \overset{R \;\;\;}{\int_{M}}(R-\bar R)dg=0.
\label{DN.7}\end{equation}
This means that we can use the result of \cite{Ji-Mazzeo-Sesum} to solve
\eqref{DN.6}.  Restricting to $M\subset \bbM$ and adding an appropriate
constant gives the desired potential function.

\end{proof}

With these potential functions, it is then easy to get long-time existence for
the normalized Ricci-flow converging to a constant scalar curvature metric
as $t\to +\infty.$

\begin{theorem} \label{thm:RicciFlow}
\begin{itemize}

\item[(i)](Ji-Mazzeo-\v{S}e\v{s}um) Suppose that $\pa_{\fun}\bar M=\emptyset$ (finite volume case) and that $\chi(M)<0.$  Let $g_{0}$ be a metric on $M$ as in Proposition~\ref{le.1}.  Then the solution $g(t)=e^{\omega(t)}g_{0}$ to
the normalized Ricci flow \eqref{RicciFlowEq} with $\sC= \bar R$ with initial metric $g_{0}$ exists
for all $t>0$ and converges exponentially fast to a complete metric
of constant negative curvature in its conformal class.
\item[(ii)]
Suppose that $\pa_{\fun}M\ne \emptyset$ and that $g_{0}$ is a metric on
$M$ satisfying \eqref{ModelFun} and
\eqref{ModelCusp}
with $\varphi\in \mathcal{C}^{4}_{\Fhc}(M)$ satisfying
\eqref{se.3b} for
some $\delta_{\fun}>2$ and $\delta_{\hc}>0.$  Then the solution $g(t)=e^{\omega(t)}g_{0}$ to
the normalized Ricci flow \eqref{RicciFlowEq} with $\sC<0$ and
with initial metric $g_{0}$ exists
for all $t>0$ and converges exponentially fast to a complete metric
of constant negative curvature in its conformal class.
\item[(iii)] If we assume that ${}^{R}\!\!\Area(g_{0})\chi(M)<0$ and that
$g_{0}$ is a smooth metric as in (ii) but with $\delta_{\fun}=2$ (instead of $\delta_{\fun}>2$) ,then the same result holds with $\sC=\bar R$.
\end{itemize}
\label{le.11}\end{theorem}
\begin{remark}
By Proposition~\ref{cbdf.1}, we know also that $\omega(\infty)\in \CI(\overline{M})$.
\end{remark}
\begin{proof}Statement (i) of the theorem is the
result of Ji-Mazzeo-\v{S}e\v{s}um \cite{Ji-Mazzeo-Sesum}.
With the potential function of Proposition~\ref{le.1}, the proof of statement (ii) is basically
the same as the one originally given by Hamilton \cite{Hamilton} in the
compact case and by \cite{Ji-Mazzeo-Sesum} in the cusp case.
We will repeat it for the convenience of the reader.

Thus, let $g_{0}$ be as in statement (ii) of the theorem and
let $g(t)$ be the solution of \eqref{RicciFlowEq} with $\sC<0.$
Let $f(t)$ denote the potential function of \eqref{le.1a}
associated to the metric $g(t).$
As in \cite{Hamilton}, one computes that
\begin{equation}
    -\Delta \frac{\pa f}{\pa t}= -\Delta(-\Delta f+ \sC f)
\label{ep.1}\end{equation}
Now, if $\psi_{i}= c_{i}\log x_{i}$ with $c_{i}:[0,T]\to \bbR$ for
$i\in\{1,\ldots, n_{\fun}\}$ are the functions such that
\[
        \tilde{f}= f- \sum_{i=1}^{n_{\fun}}\psi_{i} \in
\mathcal{C}^{2,\alpha}_{\Fhc}(M),
\]
one sees from \eqref{pot.2b} \eqref{bph.4} and \eqref{FirstLine} that the
evolution equation for $\psi_{i}$ is $\frac{\pa \psi_{i}}{\pa t}= \sC \psi_{i}$
so that from \eqref{ep.1}, we have
\begin{equation}
-\Delta \frac{\pa \tilde{f}}{\pa t}= -\Delta(-\Delta \tilde{f} +\sC \tilde{f}
-\sum_{i=1}^{n_{\fun}}\Delta \psi_{i}).
\label{ep.2}\end{equation}
Since each term satisfies the Neumann boundary conditions at
$Y_{2},\ldots,Y_{n_{\fun}}$ and the Dirichlet boundary condition at
$Y_{1}$ (modulo a constant for $\Delta \psi_{1}$),
we see by considering the corresponding equation on $\bbM_{2}$
that there exists a function $K:[0,T]\to \bbR$ such that
\begin{equation}
\frac{\pa \tilde{f}}{\pa t}= -\Delta \tilde{f} +\sC \tilde{f}
-\sum_{i=1}^{n_{\fun}}\Delta \psi_{i}+K(t), \quad
\frac{\pa f}{\pa t}= -\Delta f+ \sC f +K(t).
\label{ep.3}\end{equation}
The trick is then to consider
the function
\begin{equation}
  h= -\Delta_{g(t)}f +|\nabla f|^{2}_{g(t)}.
\label{ec.1}\end{equation}
As in \cite{Hamilton}, one computes that
\begin{equation}
    \pa_t h = - \Delta_{g(t)} h - 2|Z|^2 + \sC h\le -\Delta_{g(t)}h +\sC h,
\label{evh}\end{equation}
where $Z$ is the trace-free part of the second covariant derivative of $f.$
From the construction of the potential function $f$, the function $h(t)$ as an
asymptotic value $h_{i}(t)$ in the end $Y_{i}$ which is determined by the
function $\psi_{i}$.  This value is uniformly bounded in $t$ (as are
$r_{i}(t)$ in \eqref{FormulaAsymCurv} and $c_{i}(t)$ in \eqref{se.11}).  This
means we can therefore apply the maximum principle to \eqref{evh} to get that there exists a constant $K$ such that
$h_{\max}(t) \leq K e^{\sC t}$ for all $t$.  This implies that
\begin{equation}\label{Gronwall}
    R_t = h - |\nabla f|^2 + \sC
    \leq Ke^{\sC t} + \sC.
\end{equation}

To get a lower bound for $R-\sC,$ suppose that the minimum of
$R-\sC$ on $M$
becomes negative at a certain time $t_{0}$ (if not we get $R-\sC\ge 0$ as a lower bound).  If
$R_{\min}$ is the minimum of the curvature at this time, then
we have that $R_{\min}\le \sC.$  By \eqref{FirstLine}, the only way
the curvature can blow up is if the minimum of $R$ becomes very negative and is
attained in the interior of $M.$  In that case, from the
evolution equation of $R,$ we get that
\begin{equation}
   \frac{d}{dt}R_{\min}\ge R_{\min}(R_{\min}-\sC)\ge \sC(R_{\min}-\sC) ,
\label{ec.4}\end{equation}
from which we deduce that $R_{\min}-\sC\ge Ce^{\sC t}$ for some constant
$C.$  Combining with our upper bound, we see that
\begin{equation}
    Ce^{\sC t} \le R-\sC \le K e^{\sC t}.
\label{ec.5}\end{equation}
In particular, this gives an a priori bound on the curvature from which
we get long-time existence.  In fact, from \eqref{ec.5}, we also see that
$R$ converges to the constant $\sC$ exponentially fast as
$t\to +\infty.$  Looking at covariant derivatives of the
evolution equation of $R$ in \eqref{EvolR} and bootstrapping gives
the corresponding a priori estimates for all higher derivatives of $R.$
Integrating the flow, this gives that $g(t)$ converges
exponentially fast to a metric $g_{\infty}$ with constant scalar curvature
$\sC.$

For statement (iii), we use instead the potential function
of Proposition~\ref{DN.4} and proceed in a similar way.  We leave the details to the reader.

\end{proof}

\section{Ricci flow and the determinant of the Laplacian}\label{max.0}

Given a $\Fhc$ metric $g_0$ on a non-compact surface $M,$ we can now optimize the determinant of its Laplacian within its conformal class.
For surfaces of finite area, the analysis of the determinant on closed surfaces in \cite{OPS2} applies easily to the determinant defined with renormalized traces.
For infinite area surfaces the situation is more delicate, but also includes
situations where the Euler characteristic is nonnegative. The main difficulty is the fact that the renormalized integral of a positive density need not be positive (e.g., \S\ref{sec:RenArea}).
To deal with this, we shall have to impose restrictions on the value of the asymptotic curvature in the funnel ends of the metrics we consider.
Notice also that our definition of the determinant of the Laplacian for a $\Fhc$ metric in \S\ref{sec:RDet} requires that the factors $e^{\varphi}$ in \eqref{ModelFun} and \eqref{ModelCusp} be in $\CI(\bar M),$ since we make use of the constructions of the heat kernel in \cite{Albin2} and \cite{Vaillant}.
Thanks to Proposition~\ref{PreserveSmoothness}, this will be preserved along the flow provided
$\varphi$ is also in $\CI_{\Fhc}(M).$  We will also need to make a careful choice of the coordinates $(x_i,\theta_i)$ in \eqref{ModelFun} and \eqref{ModelCusp} adapted to the hyperbolic metric in the conformal class.

\begin{theorem}
Let $g_{\infty}$ be a hyperbolic $\Fhc$ metric on a non-compact surface $M$.  Fix the coordinates $(x_i,\theta_i)$ near each boundary component of $\overline{M}$ in such a way that \eqref{ModelFun} and \eqref{ModelCusp} hold for $g_{\infty}$ with $\varphi\equiv 0$.
Suppose $g_0$ is another $\Fhc$ metric in the conformal class of $g_{\infty}$ with totally geodesic ends on a non-compact surface $M$ satisfying \eqref{ModelFun} and \eqref{ModelCusp} with $\varphi \in \CI(\bar M) \cap \cC^{\infty}_{F,\hc}(M).$
Assume that
\begin{equation}\label{ExtraAssumption}
    {}^R\!\!\Area(g_0) = -2\pi\chi(M)\quad \mbox{and} \quad
    r_{i}=-2 \; \mbox{for}\; i\in \{1,\ldots,n_{\fun}\}.
\end{equation}
When $\chi(M)=0,$ assume that there is at least one funnel end and that $\varphi$ satisfies \eqref{se.3b} with $\delta_{\fun}>2.$
Then among all $\Fhc$ metrics $g_0$ with totally geodesic ends, conformal to $g_{\infty},$ and satisfying \eqref{ExtraAssumption} the determinant of the Laplacian is greatest at the hyperbolic metric $g_{\infty}$.
\label{max.1}\end{theorem}

\begin{proof}
Assume first that $\chi(M)\ne 0.$  Consider normalized Ricci flow starting at $g_0$ with normalization constant
\begin{equation*}
    \sC = \bar R = \frac{4\pi\chi(M)}{{}^R\!\!\Area(g_0)}<0.
\end{equation*}
We know from Lemma \ref{lem:Vol} and Theorem \ref{thm:RicciFlow} that, because of \eqref{ExtraAssumption}, this flow exists for all time, preserves the renormalized area, and converges to a hyperbolic metric.
We know from Proposition~\ref{PreserveDecay} that $\omega'(t) = \omega_0(t) + \wt\omega(t)$ with $\wt\omega(t) = \cO(x^2_{\fun})$ and
$\omega_{0}\in\CI([0,+\infty))$.  We also know from Proposition~\ref{PreserveSmoothness} that
$\omega(t)\in \CI(\bar M)\cap \mathcal{C}^{\infty}_{\Fhc}(M).$  Keeping in mind Lemma~\ref{lem:Vol},
we can therefore apply the Polyakov formula of Theorem~\ref{thm:Polyakov} to get
\begin{equation}
\begin{split}
    \pa_t \log \det(\Delta_{g(t)})
    &= -\frac1{24\pi} \Rint \omega'(t) R_t \dvol_{t} \\
    &= \frac1{24\pi} \Rint (R_t - \sC)^2 \dvol_{t} + \frac{\sC}{12\pi} \Rint (R_t - \sC) \dvol_t \\
    &= \frac1{24\pi} \int (R_t - \sC)^2 \dvol_{t} \ge 0.
\end{split}
\label{rp}\end{equation}
where the  last integral does not need renormalization because of \eqref{ExtraAssumption}.   This shows the determinant is increasing along the flow.  This will prove the theorem provided we show the determinant converges to the determinant $\Delta_{g(\infty)}$, where $g(\infty)$ is the metric of constant negative scalar curvature towards which the flow is converging.  To see this, notice that if
\[
    L(t)= \int (R_{t}-\sC)^{2}dA_{t},
\]
then using the evolution equations \eqref{EvolR}, we compute that
\begin{equation*}
\begin{aligned}
\pa_{t} L(t)&=  \int  ( 2 (R_{t}-\sC)\pa_{t}R)dA_{t} + \int (R_{t}-\sC)^{2} \pa_{t}(dA_{t})  \\
   &= \int (2(R_{t}-\sC)(-\Delta R_{t} +R_{t}(R_{t}-\sC))dA_{t} - \int (R_{t}-\sC)^{3}dA_{t} \\
   &= -2\int |\nabla (R_{t}-\sC)|^{2}dA_{t} + \sC\int(R_{t}-\sC)^{2}dA_{t} + \int(R_{t}-\sC)^{2}R_{t}dA_{t}.
\end{aligned}
\end{equation*}
The first term is obviously negative, while the third terms is negative provided $R_{t}$ is negative, which is true for $t$ sufficiently large.  Thus, for $t$ sufficiently large, we get that,
\[
      \pa_{t} L(t) \le \sC L(t),
\]
which implies $L(t)$ converges to $0$  exponentially fast as $t$ goes to infinity.  This means  that $ \det(\Delta_{g(t)})$ converges to $\det(\Delta_{g(\infty)})$, from which the theorem follows.

When $\chi(M)=0,$ consider the normalized Ricci flow starting at $g_0$ with normalization constant
\begin{equation*}
    \sC  = -2.
\end{equation*}
We know from Theorem \ref{thm:RicciFlow} that this flow exists for all time, stays within the class of metrics we are considering, and converges to a hyperbolic metric.  Clearly, \eqref{rp} still holds in this case since
\[
\overset{R}{\int}(R_{t}+2)dg_{t}= 4\pi \chi(M)+2 ({}^{R}\!\!\Area)(g_{t})=0.
\]
Since $R_t + 2 = \cO(x^{2}),$ we can prove the theorem following the same argument as in the previous case.

Finally, since the functional  $F(\omega)$ in \eqref{convex.1} is concave when $\omega$ and $\nabla \omega$ are in $L^{2},$  notice that this implies
the maximum of the determinant is unique and confirms that among all $F,\hc$ metrics $g$ conformal to
$g_{\infty}$ and satisfying the hypotheses of the theorem, there is a unique hyperbolic metric.
\end{proof}

\bibliography{RfdLncs}
\bibliographystyle{amsplain}

\end{document}